\documentclass{amsart}
\usepackage{amsmath,amssymb,bbm}
\usepackage[dvips]{epsfig}
\usepackage[all, knot]{xy}
\usepackage{hyperref}
\SelectTips{cm}{}
\xyoption{arc}
\xyoption{2cell}
\UseAllTwocells

\newcommand \uriicell [1] {\urtwocell \omit {#1}}

\newtheorem{theorem}{Theorem}
\numberwithin{theorem}{subsection}

\newtheorem{lemma}[theorem]{Lemma}
\newtheorem{varexample}[theorem]{Example}

\newtheorem{definition}[theorem]{Definition}
\newtheorem{varremark}[theorem]{Remark}
\newenvironment{remark}{\begin{varremark}\em}{\em\end{varremark}}

\newcommand{\opname}[1]{\operatorname{#1}}

\newcommand{\catname}[1]{\boldsymbol{\opname{{#1}}}}
\newcommand{\C}{\catname{C}}
\newcommand{\Cop}{\catname{C}^{\opname{op}}}
\newcommand{\Cat}{\catname{Cat}}

\newcommand{\Bicat}{\catname{Bicat}}
\newcommand{\Hor}{\catname{Hor}}
\newcommand{\Ver}{\catname{Ver}}
\newcommand{\Squ}{\catname{Squ}}
\newcommand{\iiCob}{\catname{2Cob}}
\newcommand{\nCobi}{\catname{nCob}}
\newcommand{\nCob}{\catname{nCob}_2}
\newcommand{\br}[1]{\langle {#1} \rangle}

\newcommand{\Cspan}{\opname{Span}(\catname{C})}
\newcommand{\Cosp}{\opname{Cosp}}
\newcommand{\CCosp}{\opname{Cosp}(\catname{C})}
\newcommand{\iiCCosp}{\opname{2Cosp}(\catname{C})}

\newcommand{\Obj}{\opname{Obj}}
\newcommand{\Mor}{\opname{Mor}}
\newcommand{\Ob}{\catname{Obj}}
\newcommand{\M}{\catname{Mor}}
\newcommand{\B}{\catname{2Mor}}
\newcommand{\db}{double bicategory}
\newcommand{\dbs}{double bicategories}
\newcommand{\vdb}{Verity double bicategory}
\newcommand{\vdbs}{Verity double bicategories}

\newcommand{\ra}{\mathop{\rightarrow}}

\newcommand{\rightarrowlim}{\mathop{\rightarrow}\limits}

\newcommand{\leftarrowlim}{\mathop{\leftarrow}\limits}

\begin{document}
\raggedbottom

\title{Double Bicategories and Double Cospans}

\author{Jeffrey C. Morton}

\address{Mathematics Department\\ University of Western Ontario \\ {\tt{jeffrey.c.morton@gmail.com}}}


\begin{abstract}Interest in weak cubical $n$-categories arises in
various contexts, in particular in topological field theories.  In
this paper, we describe a concept of \textit{\db} in terms of
bicategories internal to $\catname{Bicat}$.  We show that in a special
case one can reduce this to what we call a \textit{\vdb}, after
Domenic Verity.  This is a weakened version of a double category, in
the sense that composition in both horizontal and vertical directions
satisfy associativity and unit laws only up to (coherent)
isomorphisms. We describe examples in the form of double bicategories
of ``double cospans'' (or ``double spans'') in any category with
pushouts (pullbacks, respectively).  We also give a construction from
this which involves taking isomorphism classes of objects, and gives a
{\vdb} of double cospans.  Finally, we describe how to use a minor variation on this to describe cobordism of manifolds with boundary.
\end{abstract}
\maketitle

\section{Introduction}

The need to generalize the concept of a category was implicit from the
beginning of the subject. Saunders Mac Lane stated that the concept of
category was introduced to study not categories themselves, nor even
functors from one category to another, but natural transformations
between functors, which are naturally seen as 2-morphisms in a
2-category of all categories.  This was an early seed of the notion of
higher categories.  Once explicitly recognized, however, the concept
proved to be ambiguous.

There has been considerable work toward a general definition of a
(weak) $n$-category.  This has $(n+1)$ layers of structure, including
objects, morphisms between objects, 2-morphisms between morphisms, and
so on up to $n$-morphisms.  Several possible alternative
definitions exist, as discussed by Cheng and Lauda \cite{chenglauda},
and by Leinster \cite{leinster2}.  One of the features which varies
among such definitions is the \textit{shape} of higher-dimensional
morphisms, with different choices suitable to different applications.
Our aim in this paper is to develop one particular notion of higher
category, in particular a \textit{\db}, which we shall define.  We also
show that there is a broad class of examples of this type in the form
of \textit{double cospans}.

The author's original motivation here was to describe rigorously a
bicategory of \textit{cobordisms with corners}.  The most natural
development of this idea turned out to be a special case of such
double cospans.  This in turn made it clear that the most natural
structure for such things is not a bicategory, but the {\dbs}
discussed here.  However, as we will prove in Theorem \ref{thm:equiv},
given a {\db} satisfying some simple conditions, one can get a
bicategory, which is a better-understood and simpler structure.  Our
class of double span examples can be made to be of this type.  The
development of the topological material involved in cobordisms with
corners will be given in a companion paper, but here we aim to be
accessible to readers of that paper seeking background, and thus will
give a relatively expository description of {\dbs} and their double
cospan examples.

The related concept of a ``weak double category'', or ``pseudo double
category'' has also been defined (for further discussion, see e.g
Marco Grandis and Robert Par\'e \cite{GP1}, Thomas Fiore \cite{fiore},
or Richard Garner \cite{garner}).  In this setting, the weakening only
occurs in only one direction, say the horizontal.  That is, the
associativity of composition, and unit laws, in the vertical direction
apply only up to certain higher \textit{associator} and
\textit{unitor} isomorphisms.  In the horizontal direction, category
axioms hold strictly.  In fact, this must be so when weakening uses
just the square 2-cells of the double category.  For the composition
in a {\db} to be weak in both directions, it must be that the
associator isomorphisms are (globular) 2-morphisms, rather than
(square) 2-cells.

We note here that the feature that one direction is strict also
appears in the \textit{weak $n$-cubical categories} discussed by
Grandis (\cite{grandis1}, \cite{grandis3}), but that these are well 
defined for any dimension $n$.  In particular, they are defined so as to 
have one direction in which composition is strict, while all others are 
weak.  However, a {\vdb} can be taken to be a weak 3-cubical category with 
no nontrivial morphisms in the strict direction; on the other hand, a weak 
2-cubical category can be seen as a {\vdb} in which the composition in one direction happens to be strict.  For some purposes, this asymmetry is useful,
but for our motivating application to cobordisms with corners, we want to
define a ``fully'' weak cubical 2-category.  We shall comment on the structures
described by Grandis again when we consider double cospans, which give 
examples of both {\dbs} and weak $n$-cubical categories.

In Section \ref{sec:prelim} we briefly describe some of the necessary
category-theoretic background for readers who may be unfamiliar with
it.  This includes the concepts of \textit{enrichment} and
\textit{internalization} in category theory, which give rise to
bicategories and double categories respectively.  The concept of a
{\db} is a mutual generalization of that of \textit{bicategory} and of
\textit{double category}, and is best understood in this light.

Bicategories gave the first precise, explicit notion of weak higher
categories.  They were described by B\'enabou
\cite{benabou} in 1967, introducing the concept of 2-morphisms between
morphisms:
\begin{equation}\label{xy:2morphism}
  \xymatrix{
    x \ar@/^1pc/[r]^{f}="0" \ar@/_1pc/[r]_{g}="1" & y \\
      \ar@{=>}"0"+<0ex,-2ex> ;"1"+<0ex,+2ex>^{\alpha}
  }
\end{equation} Equations in the axioms for a category are replaced
by 2-isomorphisms, which themselves satisfy coherence laws given in
equations.  This is known as \textit{weakening}.

The original example used to illustrate this concept was the
bicategory of \textit{spans} in a suitable category $\catname{C}$,
namely diagrams of the form:
\begin{equation}
X \rightarrow S \leftarrow Y
\end{equation}
There is a natural concept of a map of spans, and an operation of
composition for spans which is not strictly associative---rather, it
is only associative up to isomorphism.  Weakening thus appeared
naturally in the setting of spans.  The double cospans defined here
lead to weakening in just the same way, but the concept being weakened
is that of a double category.

Double categories, introduced by Ehresmann \cite{ehresmann}, have
objects, horizontal and vertical morphisms which can be represented
diagrammatically as edges, and squares:
\begin{equation}\label{xy:square}
 \xymatrix{
  x \ar[r]^{\phi} \ar[d]_{f}  & x' \ar[d]^{f'} \\
  y \uriicell{F} \ar[r]_{\hat{\phi}} & y'
 }
\end{equation}
These can be composed in geometrically obvious ways.

Double categories distinguish between horizontal and vertical
1-morphisms, which in general can only be composed with other
morphisms of the same type.  On the other hand, 2-morphisms are
``squares'', with both horizontal and vertical source and target,
which can be composed in either direction with other squares having a
common boundary. 

Moskaliuk and Vlassov \cite{dblcatmp} discuss the application of
double categories to mathematical physics, and particularly to
topological quantum field theories (TQFT's), and to dynamical systems
with changing boundary conditions---that is, with inputs and outputs.
Kerler and Lyubashenko \cite{KL} describe ``extended'' TQFT's as
``double pseudofunctors'' between double categories.  This formulation
involves, among other things, a double category of cobordisms with
corners.  This sort of topological category has manifolds for objects,
and manifolds with boundary or with corners as higher morphisms.  This
makes it possible to describe systems with changing boundary
conditions, and the most natural way to do this is by allowing both
initial and final states, and changing boundary conditions, as part of
the boundary in a more general sense.  This is one of the main
motivations for the concepts we describe here, and we shall return to
it in a subsequent paper.  Double categories are too strict to be
really natural for our purpose, however.  Composition in a double
category must be strictly associative, and in order to achieve this,
one considers only \textit{equivalence classes} of cobordisms as
morphisms.

Thus, the principle here is to weaken the definition of a double
category.  This had previously been done in the definition of a
\textit{pseudocategory}, as described, for instance, by Fiore
\cite{fiore}.  However, in a pseudocategory, just one direction of
composition is weak: that is, the associative and unit laws satisfied
by composition are replaced by associator and unitor isomorphisms.  In
a {\db}, composition is weak in both directions.

In Section \ref{sec:doublebicat} we introduce {\dbs} using a form of internalization, analogous to that which gives double categories as categories internal to $\Cat$.  In Section \ref{sec:vdb}, we
describe a somewhat different concept of {\db}, due to Dominic Verity
(which we denote a \textit{\vdb} for clarity), a structure which also
has both horizontal and vertical bicategories, and square 2-cells,
with weak composition in both directions.  In Section
\ref{sec:decatfy}, we explain how a special case of our {\dbs} can be
reduced to Verity's definition.  In turn, we show in Section
\ref{sec:equiv} how a {\vdb} satisfying certain conditions, in turn
yields a bicategory in the usual sense.  Thus, this presents a series
of increasingly manageable simplifications.

We finish in Section \ref{sec:dblspan} by describing a rather broad
general class of examples of {\dbs}, which arise in rather the same
way as the fact that $\Cspan$ is a bicategory.  A ``double cospan'' in
a category with pushouts is a diagram of the following form:
\begin{equation}
  \xymatrix{
    X_1  \ar[r] \ar[d] & S \ar[d] & X_2 \ar[l] \ar[d] \\
    T_1 \ar[r]  & M  & T_2   \ar[l] \\
    X'_1 \ar[u] \ar[r] & S' \ar[u] & X'_2 \ar[u] \ar[l]
  }
\end{equation}These diagrams can be composed horizontally and
vertically, and in either case composition is by pushout, just as with
ordinary spans or cospans.  In Section \ref{sec:dcexample} we describe
double cospans and their composition in detail, and show that they
naturally form a {\db} in our original sense.  In Section
\ref{sec:vdbcosp} we show how reducing to certain natural equivalence
classes of double cospans yields a {\vdb}, following the procedure in
Section \ref{sec:decatfy}.

\section{Bicategories and Double Categories}\label{sec:prelim}

Since this paper is intended as a companion to another of a more
topological nature, we will recall for the reader with less
category-theoretic background some relevant ideas about bicategories
and double categories.  Other readers may wish to skip to Section
\ref{sec:doublebicat} when we introduce double bicategories.

We want to weaken the notion of a double category.  Weakening a
concept $X$ in category theory generally involves creating a new
concept in which defining equations in the original concept (such as
associativity) are replaced by specified isomorphisms.  Thus, one says
that the defining equations hold with equality in a \textit{strict
$X$} and hold only ``up to'' isomorphism in a \textit{weak $X$}.

Before describing our weakened concept of double category, we recall
how this process works, and examine the strict form of the concept we
want to weaken.  So we begin by recalling some facts about
bicategories, to illustrate weakening, and double categories, to
provide a starting point.

\subsection{Bicategories}\label{sec:bicat}

A \textbf{bicategory} is a ``weak globular 2-category''.  That is, if
$\catname{B}$ is a bicategory, and $x,y \in \catname{B}$, then
$\hom(x,y) \in \Cat$, allowing isomorphisms between morphisms where
formerly we had equations.  The morphisms in $\hom(x,y)$ are thought
of as ``2-morphisms'' in $\catname{B}$.  Moreover, the strict version
of a bicategory, usually called a ``2-category'', has the same unit
and associativity axioms as a category.  However, the weak form
replaces these with 2-isomorphisms satisfying some coherence
properties.  So in particular, we have the following definition, due
to B\'enabou \cite{benabou}, and discussed in more detail, for
instance, in \cite{leinster}.

\begin{definition} A \textbf{bicategory} $\mathcal{B}$ consists of the
following data:
\begin{itemize}
\item A collection of \textbf{objects} $\Ob$
\item For each pair $x,y \in \Ob$, a \textit{category} $\hom(x,y)$
      whose objects are called \textbf{morphisms} of $\mathcal{B}$
      and whose morphisms are called \textbf{2-morphisms} of
      $\mathcal{B}$
\item For each object $x \in \Ob$, an \textit{identity} $1_x \in
      \hom(x,x)$
\item For each triple $x$, $y$, $z$ of objects, a \textbf{composition}
      functor $\circ : \hom(x,y) \times \hom(y,z) \rightarrow
      \hom(x,z)$
\item For each composable triple $f$, $g$, $h$ of morphisms, a
      2-isomorphism (i.e. invertible 2-morphism) $\alpha_{f,g,h} : h
      \circ (g \circ f) \rightarrow (h \circ g) \circ f$ called the
      \textbf{associator}
\item For each morphism $f: x \rightarrow y$, left and right
      \textbf{unitor} 2-isomorphisms $l_f : 1_y \circ f \rightarrow f$
      and $r_f : f \circ 1_x \rightarrow f$
\end{itemize}

The associator is subject to the Pentagon identity, namely that the
following diagram commutes for any 4-tuple of composable morphisms
$(f,g,h,j)$:
\begin{equation}\label{eq:pentagonid}
\xy
 (0,20)*+{(f \circ g) \circ (h \circ j)}="1";
 (35,4)*+{f \circ (g \circ (h \circ j))}="2";
 (23,-20)*+{ \quad f \circ ((g \circ h) \circ j)}="3";
 (-23,-20)*+{(f \circ (g \circ h)) \circ j}="4";
 (-35,4)*+{((f \circ g) \circ h) \circ j}="5";
     {\ar^{a_{f,g,h \circ j}}     "1";"2"}
     {\ar_{1_f \circ a _{g,h,j}}  "3";"2"}
     {\ar_{a _{f,g \circ h,j}}    "4";"3"}
     {\ar_{a _{f,g,h} \circ 1_j}  "5";"4"}
     {\ar^{a _{f \circ g,h,j}}    "5";"1"}
\endxy
\\
\end{equation}
Also, the unitors and associator make the following commute for all
composable $g,f$:
\begin{equation}\label{eq:unitorlaws}
  \xymatrix{
    (g \circ 1_y) \circ f \ar[r]^{a_{g,1_y,f}} \ar[d]_{r_g \circ 1_f} & g \circ (1 \circ f) \ar[dl]^{1_g \circ l_f}\\
    g \circ f & \\
  }
\end{equation} (where $y = t(f) = s(g)$).
\end{definition}

\begin{remark} This is a compact definition of a bicategory, but it is
possible to describe the same data in different ways, which will be
more directly relevant to subsequent discussion of {\dbs}.  In
particular, this definition is related to the definition of a
\textit{strict} bicategory (a \textit{2-category}) as a
\textit{category enriched in} $\catname{Cat}$.  That is, for any
objects $x$ and $y$, there is a category $\hom(x,y)$.  However, there
is a more elementary, although perhaps less elegant, way of describing
the data of a bicategory.

One can form the collection $\M = \coprod \opname{ob}(\hom(x,y))$ of
all morphisms of $\mathcal{B}$, and $\B = \coprod
\opname{mor}(\hom(x,y))$ of all 2-morphisms of $\mathcal{B}$.  Then
the existence of composition functors imply that there is, just as in
categories, a partially defined composition function $\circ : \M \times
\M \rightarrow \M$, which is defined for pairs $(f,g)$ for which
$t(f)=s(g)$ (and two such functions giving composition of
2-morphisms).  These functions will have properties determined by the
fact that they must give composition functors as defined above.  The
existence of identity morphisms means that there is an
\textbf{identity} map $i :\Ob \rightarrow \M$, and this satisfies the
usual relations with the source and target maps, and the composition.

As well as the source and target maps
\[
s,t: \M \rightarrow \Ob
\]
given in this definition, there are the source and target maps in each
$\hom(x,y)$.  These imply the existence of $s,t : \B \rightarrow \M$,
with the property that for any 2-morphism $\alpha$, $s(s(\alpha)) =
s(t(\alpha))$, and $t(s(\alpha))=t(t(\alpha))$.  A similar condition
will apply to double categories, and indeed {\dbs}, as we shall see.
Together these describe the picture summarized in the diagram
(\ref{xy:2morphism}), depicting 2-morphisms as 2-dimensional cells
between arrow-shaped 1-morphisms.
\end{remark}

Jean B\'enabou \cite{benabou} introduced bicategories in a 1967 paper,
and one broad class of examples introduced there comes from the notion
of a \textit{span}.  Since we will want to use a similar construction
later, we remark on this here:

\begin{definition} (\textbf{B\'enabou}) Given any category $\C$, a \textbf{span}
$(S,\pi_1,\pi_2)$ between objects $X_1,X_2 \in \C$ is a diagram in
$\C$ of the form
\begin{equation}
  \xymatrix{
    P_1 & S \ar[l]_{\pi_1} \ar[r]^{\pi_2}  & P_2
  }
\end{equation} Given two spans $(S,s,t)$ and $(S',s',t')$ between
$X_1$ and $X_2$ a \textbf{morphism of spans} is a morphism $g:S
\rightarrow S'$ making the following diagram commute:
\begin{equation}\label{eq:spanmorph}
  \xymatrix{
        & S \ar[dl]_{\pi_1} \ar[dr]^{\pi_2} \ar[d]^{g} &   \\
    X_1 & S' \ar[l]^{\pi'_1} \ar[r]_{\pi'_2} & X_2 \\
  }
\end{equation}

Composition of spans $S$ from $X_1$ to $X_2$ and $S'$
from $X_2$ to $X_3$ is given by pullback: that is, an
object $R$ with maps $f_1$ and $f_2$ making the following diagram
commute:
\begin{equation}\label{eq:spanpullback}
  \xymatrix{
      &   & R \ar[dl]_{f_1} \ar[dr]^{f_2} &  & \\
      & S \ar[dl]_{\pi_1} \ar[dr]^{\pi_2} &   & S' \ar[dl]_{\pi'_2} \ar[dr]^{\pi'_3} & \\
    X_1 &  & X_2 &   & X_3 \\
  }
\end{equation} which is terminal among all such objects.  That is,
given any other $Q$ with maps $g_1$ and $g_2$ which make the analogous
diagram commute, these maps factor through a unique map $Q \rightarrow
R$.  $R$ becomes a span from $X_1$ to $X_3$ with the maps $\pi_1 \circ f_1$
and $\pi_2 \circ f_2$.
\end{definition}

The span construction has a dual concept:

\begin{definition} A \textbf{cospan} in $\C$ is a span in $\Cop$,
morphisms of cospans are morphisms of spans in $\Cop$, and composition
of cospans is given by pullback in $\Cop$---that is, by pushout in
$\C$.
\end{definition}

\begin{remark}\label{rem:spanbicat} (\textbf{\cite[ex.\ 2.6]{benabou}})
\label{thm:spanbicat}Given any category $\C$ with all limits, there is
a bicategory $\Cspan$, whose objects are the objects of $\C$, whose
$hom$-sets of morphisms $\Cspan(X_1,X_2)$ consist of all spans between
$X_1$ and $X_2$ with composition as defined above, and whose
2-morphisms are morphisms of spans.  $\Cspan$ as defined above forms a
bicategory (dually, there is a bicategory $\Cosp(\C)$ of cospans).

One should note that there is some choice in the precise definition of
this bicategory since pushout (or pullback) is only defined up to
isomorphism.  However, one can make a particular choice of pushout (or
pullback) as a given composite, and given this choice get
corresponding associators and unitors.  Different choices of composite
will of course give different such maps.  However, all such choices
are equivalent. This is due in part to the fact that the pullback is a
universal construction (universal properties of $\Cspan$ are discussed
by Dawson, Par\'e and Pronk \cite{unispan}).

We briefly describe the proof of B\'enabou that $\Cspan$ is a
bicategory:

The identity for $X$ is $X \leftarrowlim^{id} X \rightarrowlim^{id}
X$, which has an obvious unitor whose properties are easy to check.

The associator arises from the fact that the pullback is a
\textit{universal} construction.  Given morphisms $f: X \rightarrow
Y$, $g: Y \rightarrow Z$, $h: Z \rightarrow W$ in $\Cspan$, the
composites $((f \circ g) \circ h)$ and $(f \circ (g \circ h))$ are
pullbacks consisting of objects $O_1$ and $O_2$ with maps into $X$ and
$W$.  The universal property of pullbacks gives an isomorphism between
$O_1$ and $O_2$ as follows.

The universal property of pullback means that any object with maps
into the objects $f$ and $(g \circ h)$ will have a map into $O_2$ which they
factor through.  We have maps into the objects $f$, $g$, and $h$ from
$O_1$, and therefore a unique compatible map into $g \circ h$ by the
universal property for that pullback.  Therefore, there is again a
unique compatible map into $O_2$.  This we take to be the associator.
We notice that in particular, the same argument works in reverse, and
so the two maps we get are inverses, hence isomorphisms.

These associators satisfy the pentagon identity
since they are unique (in particular, both sides of the pentagon give
the same isomorphism).

It is easy to check that $\hom(X_1,X_2)$ is a category, since it
inherits all the usual properties from $\C$.
\end{remark}

We will generalize the construction of bicategories of spans to give
examples of double bicategories in Section \ref{sec:dblspan}.  This
development of bicategories illustrates the sort of weakening we want
to apply to the concept of a double category.  So we will first
describe the strict notion in Section \ref{sec:doublecat}, before
considering how to weaken it, in Section \ref{sec:vdb}.

\subsection{Double Categories}\label{sec:doublecat}

The concept of a double category extends that of a category in a
different way than the concept of bicategory.  Both, however, can be
visualized as having both ``arrow-like'' morphisms, and also
two-dimensional cells thought of as higher morphisms.  

Just as with bicategories, we recall the definition here first by
giving an abstract definition, then showing an equivalent, more
concrete, version.  The first definition of a bicategory highlighted
its relation to the idea of an enriched category.  Here we begin by
illustrating how double categories illustrate
\textit{internalization}, which we will return to in Section
\ref{sec:doublebicat} when describing {\dbs}.

\begin{definition}A double category is a category internal to $\Cat$.
\end{definition}

This is a generalization of the more broadly familiar terminology in
which, for instance, a group internal to $\catname{Top}$ is called a
\textit{group object} in $\catname{Top}$, or topological group.
However, not all structures we might want to internalize are
determined by a single object.  In particular, a category (by default,
internal to $\catname{Set}$) consists of not one but two sets, namely
the sets of \textit{objects} and \textit{morphisms}.  A category
internal to $\catname{C}$ (or ``in $\catname{C}$'') has two objects of
$\catname{C}$ playing the same roles.

So a category in $\Cat$ is a structure having a category $\catname{Ob}$
of objects and a category $\catname{Mor}$ of morphisms, with
functors such as $s$ and $t$ satisfying the usual category axioms.  Note 
that these axioms give conditions at both the object and morphism level,
in addition to those which follow from the fact that they are functors.

We thus have sets of objects and morphisms in $\catname{Ob}$, which
satisfy the usual axioms for a category.  The same is true for
$\catname{Mor}$.  In addition, the category axioms for the double category
are imposed on the composition and identity functors, and these must be compatible with the category axioms in the other direction.  Thus we can 
think of both the objects in $\M$ and the morphisms in $\Ob$ as morphisms
between the objects in $\Ob$.  A double category is often thought of as
including the morphisms of two (potentially) different categories on
the same collection of objects.  These are the \textit{horizontal} and
\textit{vertical} morphisms.

Here, the objects in the diagram can be thought of as objects in
$\Ob$, the vertical morphisms $f$ and $f'$ can be thought of as
morphisms in $\Ob$ and the horizontal morphisms $\phi$ and
$\hat{\phi}$ as objects in $\M$.  Vertical composition is given by
composition in $\M$, and horizontal composition by the morphism map of
the composition functor $\circ$.  (In fact, we can adopt either
convention for distinguishing horizontal and vertical morphisms).
However, we also have morphisms in $\M$.  We represent these as
2-cells, or {\textit{squares}}, like the 2-cell $S$ represented in
(\ref{xy:doublecat2cell}).  The fact that the composition map $\circ$
is a functor means that horizontal and vertical composition of square
2-cells commutes.

A double category can therefore be seen more directly.  It consists
of:
\begin{itemize}
\item a set $O$ of objects
\item \textit{horizontal} and \textit{vertical} categories, whose sets
      of objects are both $O$
\item for any diagram of the form
     \begin{equation}
      \xymatrix{
       x \ar[r]^{\phi} \ar[d]_{f} & x' \ar[d]^{f'} \\
       y \ar[r]_{\phi'} & y'
      }
\end{equation}
      a collection of \textit{2-cells}, having horizontal
      source and target $f$ and $f'$, and vertical source
      and target $\phi$ and $\phi'$
\end{itemize} along with additional data such as the source and target
maps, identities, and so forth, all satisfying category-like axioms in
both horizontal and vertical directions.  In particular, the 2-cells
can be composed either horizontally or vertically in the obvious way.
We denote a 2-cell filling the above diagram like this:
\begin{equation}\label{xy:doublecat2cell}
 \xymatrix{
  x \ar[r]^{\phi} \ar[d]_{f} & x' \ar[d]^{f'} \\
  y \ar[r]_{\phi'} \uriicell{S} & y'
 }
\end{equation} and think of the composition of 2-cells as pasting them
along an edge.  The resulting 2-cell fills a square whose boundaries
are the corresponding composites of the morphisms along its edges.

Next, in Section \ref{sec:doublebicat} we take our descriptions of
double categories and bicategories, and see how to find some common
generalizations of both.

\section{Double Bicategories}\label{sec:doublebicat}

We wish to describe a structure which is sufficient to reproduce the
various types of composition found in a double category, but in such a
way that all are weakened.  This means we should have horizontal
composition for horizontal morphisms and vertical composition for
vertical morphisms.  Square 2-cells should be composable in both directions.
Composition of morphisms in each direction is to be ``weak'', in the
sense of having associator and unitor isomorphisms rather than
associativity and unit laws.  This means there will also be
2-morphisms of some appropriate shape to act as unitors and
associators (and, of course, there will in general be other
2-morphisms as well).  In particular, in place of the mere categories
found in a double category, we have horizontal and vertical
\textit{bicategories}, with their (globular) 2-morphisms, as well as
(square) 2-cells.

The natural choice of name for such a structure is a \textit{\db}.
This term seems to have been originally introduced by Dominic Verity
\cite{verity}.  There is some ambiguity here.  By analogy with
``double category'', the term ``{\db}'' might be expected to describe
is an internal bicategory in $\Bicat$, the category of all
bicategories .  Indeed, it is what we will mean by a {\db} here.
However, this is not quite the concept given by Verity.  The two turn
out to be closely related, and both will be important, so we
will refer to {\dbs} in the sense of Verity by the term
\textit{{\vdbs}}, while reserving \textit{\db} for internal
bicategories in $\Bicat$.  For more discussion of the relation between
these, see Section \ref{sec:vdb}.

\subsection{Double Bicategories and Internalization}\label{sec:internal}

Here we present a more precise definition of the concept of a {\db} as a
bicategory internal to $\Bicat$.  This will be somewhat more complex than the analogous process for a double category, but runs along similar lines.

Thus, we will have bicategories $\Obj$, $\opname{Mor}$ and
$\opname{2Mor}$.  Then one can consider a bicategory internal to
$\Bicat$.  It is straightforward to treat $F(\Obj)$ as a horizontal
bicategory, and the objects of $\Obj$, $\opname{Mor}$ and
$\opname{2Mor}$ as forming a vertical bicategory.  Note, however, that
a diagrammatic representation of, for instance, 2-morphisms in
$\opname{2Mor}$ would require a 4-dimensional diagram element.  The
comparison can be seen by contrasting tables \ref{table:doublecatdata}
and \ref{table:doublebicatdata} in Section \ref{sec:decatfy}.

\begin{definition}
A {\db} consists of:
\begin{itemize}
\item bicategories $\Ob$ of \textbf{objects}, $\M$ of
      \textbf{morphisms}, $\B$ of \textbf{2-morphisms}
\item \textbf{source} and \textbf{target} 2-functors
  \begin{itemize}
  \item $s,t: \M \ra \Ob$
  \item $s,t: \B \ra \Ob$
  \item $s,t: \B \ra \M$
  \end{itemize}
\item  Composition 2-functors:
  \begin{itemize}
  \item $\circ : \opname{MPairs} \rightarrow \M$
  \item $\circ : \opname{HPairs} \rightarrow \B$
  \item $\cdot : \opname{VPairs} \rightarrow \B$
  \end{itemize}  satisfying the interchange law,
  where
  \begin{itemize}
  \item $\opname{MPairs} = \M \times_{\Ob} \M$
  \item $\opname{HPairs} = \B \times_{\M} \B$
  \item $\opname{VPairs} = \B \times_{\Ob} \B$ 
  \end{itemize}

  are (strict) pullbacks
\item an \textbf{associator} 2-functor
  \begin{itemize}
  \item $a : \opname{Triples} \rightarrow \opname{2Mor}$
  \end{itemize}
  where
 \begin{itemize}
  \item $\opname{Triples} = \M \times_{\Ob} \M \times_{\Ob} \M$
 \end{itemize}
\item \textbf{unitors}
  \begin{itemize}
   \item $l,r : \Ob \rightarrow \M$
  \end{itemize}
\end{itemize}
such that $a$ makes the following diagram commute:
  \[
    \xymatrix{
     \opname{Pairs} \ar[d]^{\circ} &  \opname{Triples} \ar[r]^{1 \times \circ} \ar[d]^{a} \ar[l]_{\circ \times 1}  & \opname{Pairs} \ar[d]^{\circ} \\
     \opname{Mor} & \opname{2Mor} \ar[r]_{t} \ar[l]^{s} & \opname{Mor} \\
    }
  \] and additional diagrams with the interpretation that $a$ gives
 \textit{invertible} 2-morphisms. The unitors must satisfy $s(l(x)) =
 t(l(x)) = x$ and $s(r(x)) = t(r(x)) = x$, and the associator should
 satisfy the pentagon identity (\ref{eq:pentagonid}), and the unitors should satisfy the unitor laws
(\ref{eq:unitorlaws}).
\end{definition}

We interpret these morphisms involving pullbacks (the fibred products)
as giving \textit{partially defined} composition 2-functors $\circ :
\M^2 \ra \M$, $\circ : \B^2 \rightarrow \B$ and $\cdot : \B^2
\rightarrow \B$, and associator 2-functor $a : \M^3 \rightarrow \B$.

\begin{remark}The Pentagon identity is shown in (\ref{eq:pentagonid})
for a bicategory (i.e. a bicategory in $\catname{Sets}$).  In
$\Bicat$, this holds for objects, morphisms, and 2-morphisms.  We can
express this condition formally, in any category (with pullbacks),
building from composable quadruples, so that the pentagon identity is
expressed in a commuting diagram which includes the one built by
pasting the two following diagrams together along the outside edges:
\begin{equation}\label{eq:pentagonth3}
\xy
(-30,40)*+{\opname{Mor}}="a1";
(30,40)*+{\opname{Mor}}="a4";
(0,25)*+{\opname{4tuples}}="b";
(-25,10)*+{\opname{Triples} \times \opname{Mor}}="c1";
(0,10)*+{\opname{Triples}}="c2";
(25,10)*+{\opname{Mor} \times \opname{Triples}}="c3";
(-25,-5)*+{\opname{2Mor}^2}="d1";
(0,-5)*+{\opname{2Mor}}="d2";
(25,-5)*+{\opname{2Mor}^2}="d3";
(-20,-15)*+{\opname{2Mor}}="e1";
(0,-20)*+{\opname{VTriples}}="e2";
(20,-15)*+{\opname{2Mor}}="e3";
(0,-35)*+{\opname{2Mor}}="f";
{\ar^{\pi_1\circ i}  "b";"a1"};
{\ar_{\pi_4\circ i}"b";"a4"};
{\ar@/_1pc/ "b";"c1"};
{\ar^{(id \times \circ \times id)i} "b";"c2"};
{\ar@/^1pc/ "b";"c3"};
{\ar^{a \times id} "c1";"d1"};
{\ar^{a} "c2";"d2"};
{\ar^{id \times a} "c3";"d3"};
{\ar^{\circ} "d1";"e1"};
{\ar^{p_2} "e2";"d2"};
{\ar^{\circ} "d3";"e3"};
{\ar^{p_1} "e2";"e1"};
{\ar_{p_3} "e2";"e3"};
{\ar^{\circ} "e2";"f"};
{\ar@/_6pc/_{s} "a1";"f"};
{\ar@/^6pc/^{t} "a4";"f"};
\endxy
\end{equation}
and
\begin{equation}\label{eq:pentagonth2}
\xy
(-30,40)*+{\opname{Mor}}="a1";
(30,40)*+{\opname{Mor}}="a4";
(0,25)*+{\opname{4tuples}}="b";
(-15,15)*+{\opname{Triples}}="c1";
(15,15)*+{\opname{Triples}}="c3";
(-25,0)*+{\opname{2Mor}^2}="d1";
(25,0)*+{\opname{2Mor}^2}="d3";
(-20,-15)*+{\opname{2Mor}}="e1";
(0,-20)*+{\opname{VPairs}}="e2";
(20,-15)*+{\opname{2Mor}}="e3";
(0,-35)*+{\opname{2Mor}}="f";
{\ar^{\pi_1\circ i}  "b";"a1"};
{\ar_{\pi_4\circ i}"b";"a4"};
{\ar "b";"c1"};
{\ar "b";"c3"};
{\ar^{a} "c1";"d1"};
{\ar^{a} "c3";"d3"};
{\ar^{\circ} "d1";"e1"};
{\ar^{\circ} "d3";"e3"};
{\ar^{\pi_1 \circ i} "e2";"e1"};
{\ar_{\pi_2 \circ i} "e2";"e3"};
{\ar^{\circ} "e2";"f"};
{\ar@/_6pc/_{s} "a1";"f"};
{\ar@/^6pc/^{t} "a4";"f"};
\endxy
\end{equation}

Note that diagram (\ref{eq:pentagonth3}) denotes the three 2-morphism
sequence in the pentagon, and (\ref{eq:pentagonth2}) the sequence of
two 2-morphisms.

Similar remarks apply to give ``element-free'' versions of the
interchange laws for composition of 2-morphisms and the unitor laws
shown in (\ref{eq:unitorlaws}).
\end{remark}

To fully expand this definition without assuming the concept of a
bicategory would be much longer than the form given here.  One would
have to specify nine types of data - objects, morphisms, and
2-morphisms in each of $\Ob$, $\M$, and $\B$, and describe all the
axioms in detail, such as the conditions implied by the fact that
$\circ$ and $\cdot$ are functors.  This is a rather complicated
structure, as we see in more detail when we return to it in Section
\ref{sec:decatfy} (and in particular we show the types of data implied
by this definition in Table \ref{table:doublebicatdata}).

In particular, the most natural geometric representation of a
2-morphism in $\B$ is as a 4-dimensional object.  We had hoped for a
common generalization of double categories and bicategories, each of
which is represented graphically with morphisms being cells having
dimension at most 2.  One could hope that such a structure would also
have at most 2-dimensional morphisms.  The definition of a {\vdb}
satisfies this, as we describe in Section \ref{sec:vdb}, and in
Section \ref{sec:decatfy} we show how it is related to the definition
we have given here.

First, we briefly remark that one can cast this description
of internal bicategories in terms of models of the \textit{finite
limit theory} of bicategories, $\catname{Th(\Bicat)}$.  This is a
category with finite limits, which can be described in terms of its
generators and relations.  It is generated by objects $O$, $M$ and
$B$, together with morphisms, and subject relations, as given in the
definition (where in that case these are in $\Bicat$.  A model of such
a theory in a category $\C$ with finite limits is a functor $F :
\catname{Th(\Bicat)} \rightarrow \C$.

Here we are considering \textit{strict} models of the theory of categories
in $\Cat$, and bicategories in $\Bicat$, rather than a weak model, which 
one might also consider.  In particular, $\Bicat$ is a tricategory (defined 
by Gordon, Power and Street \cite{gps}): it has objects which are
bicategories, morphisms which are (weak) 2-functors between
bicategories, 2-morphisms which are natural transformations between
bifunctors, and 3-morphisms which are modifications of such
transformations.  However, for both double categories and double
bicategories, we ignore the higher morphisms in this setting and think
of $\Bicat$ as a mere category, taking equivalence classes of
morphisms between bicategories (that is, 2-functors) where needed.  Thus,
equations in the theory are mapped to equations (not isomorphisms) in
$\Bicat$.

Even such strict models, however, are fairly complex structures, so we
now consider one way to simplify them.

\subsection{Verity Double Bicategories}\label{sec:vdb}

The following definition of a ``{\db}'' is due to Dominic Verity
\cite{verity}, and will henceforth be referred to as a {\vdb}.  It is
readily seen as a natural weakening of the definition of a double
category.  Just as the concept of \textit{bicategory} is weaker than
that of \textit{2-category} by weakening the associative and unit
laws, {\vdbs} will be weaker than double categories.

\begin{definition}\label{def:doublebicat} (\textbf{Verity}) A
\textbf{\vdb} $\catname{C}$ is a structure $\mathcal{V}$ consisting of
the following data:

\begin{itemize}
  \item a class of \textbf{objects} $\Obj$,

  \item \textbf{horizontal} and \textbf{vertical bicategories} $\Hor$
        and $\Ver$ having $\Obj$ as their objects

  \item for every square of horizontal and vertical morphisms of the form
  \begin{equation}
    \xymatrix{
      a \ar[r]^{h} \ar[d]_{v} & b \ar[d]^{v'} \\
      c \ar[r]^{h'} & d \\
    }
  \end{equation}
  a class of \textbf{squares} $\Squ$, with maps $s_h, t_h : \Squ
  \rightarrow \Mor(\Hor)$ and $s_v, t_v : \Squ \rightarrow \Mor(\Ver)$,
  satisfying an equation for each corner, namely:
  \begin{eqnarray}\label{eq:squarestmaps}
    s(s_h)&=&s(s_v) \\
    \nonumber t(s_h)&=&s(t_v) \\
    \nonumber s(t_h)&=&t(s_v) \\
    \nonumber t(t_h)&=&t(t_v)
  \end{eqnarray}
\end{itemize}
  The squares should have horizontal and vertical composition
  operations, defining the vertical composite $F \otimes_V G$
  \begin{equation}\label{eq:squarevertcomp}
    \xymatrix{
      x \ar[r] \ar[d] & x' \ar[d] \\
      y \ar[r] \ar[d] \uriicell{F} & y'\ar[d] \\
      z \ar[r] \uriicell{G} & z' \\
    } \qquad = \qquad
    \xymatrix{
      x \ar[r] \ar[d] & x' \ar[d] \\
      z \ar[r] \uriicell{F\otimes_V G} & z'\\
    }
  \end{equation}
  and horizontal composite $F \otimes_H G$:
  \begin{equation}\label{eq:squarehorizcomp}
    \xymatrix{
      x \ar[r] \ar[d] & y \ar[d] \ar[r] & z \ar[d] \\
      x' \ar[r] \uriicell{F} & y' \ar[r] \uriicell{G} & z'
    } \qquad = \qquad
    \xymatrix{
      x \ar[r] \ar[d] & z \ar[d] \\  
      x' \ar[r] \uriicell{F \otimes_H G} & z' \\
    }
  \end{equation} The composites have the usual relation to source and target
  maps, satisfy the interchange law
  \begin{equation}\label{eq:squareinterchangelaw}
    (F \otimes_V F') \otimes_H (G \otimes_V G') = (F \otimes_H G) \otimes_V (F' \otimes_H G')
  \end{equation} and there is a unit for composition of squares: 
  \begin{equation}
    \xymatrix{
      x \ar[r]^{1_x} \ar[d]_{f} & x \ar[d]^{f} \\
      y \ar[r]^{1_y} \uriicell{1_f} & y \\
    }
  \end{equation}
  (and similarly for vertical composition).

  There is a left and right action by the horizontal
  and vertical 2-morphisms on $\Squ$, giving $F \star_V \alpha$,
  \begin{equation}
    \xymatrix{
      x \ar[r] \ar[d] & y \ar[d]^{ }="1" \ar@/^2pc/[d]^{ }="0" \\
      x' \ar[r] \uriicell{F} & y' \\
      \ar@{<=}"0" ;"1"^{\alpha}
    } \qquad = \qquad
    \xymatrix{
      x \ar [r] \ar[d] & y \ar[d] \\
      x' \ar[r] \uriicell{F \star_V \alpha} & y'
    }
  \end{equation}
  (and similarly on the left) and $F \star_H \alpha$,
  \begin{equation}
    \xymatrix{
      x \ar[r]^{ }="1" \ar@/^2pc/[r]^{ }="0" \ar[d] & y \ar[d]  \\
      x' \ar[r] \uriicell{F} & y' \\
      \ar@{=>}"0" ;"1"^{\alpha}
    } \qquad = \qquad
    \xymatrix{
      x \ar [r] \ar[d] & y \ar[d] \\
      x' \ar[r] \uriicell{\alpha \star_H F} & y'
    }
  \end{equation} The actions also satisfy interchange laws:
  \begin{equation}\label{eq:actioninterchange}
    (F \otimes_H F') \star_H (\alpha \otimes_V \alpha') = (F \star_H \alpha) \otimes_h (F' \star_H \alpha')
  \end{equation}
(and similarly for the vertical case) and are compatible with composition: 
  \begin{equation}\label{eq:actioncompat}
    (F \otimes_H G) \star_V \alpha = F \otimes_H (G \star_V \alpha)
  \end{equation} (and analogously for vertical composition).  They
  also satisfy additional compatibility conditions: the left and right
  actions of both vertical and horizontal 2-morphisms satisfy the
  ``associativity'' properties,
 \begin{equation}
    \alpha \star (S \star \beta) = (\alpha \star S) \star \beta
  \end{equation}
  for both $\star_H$ and $\star_V$.  Moreover, horizontal
  and vertical actions are independent:
  \begin{equation}\label{eq:actionindep}
    \alpha \star_H (\beta \star_V S) = \beta \star_V (\alpha \star_H S)
  \end{equation}
  and similarly for the right action. 

Finally, the composition of squares agrees with the associators for
composition by the action in the sense that given three composable
squares $F$, $G$, and $H$:
\begin{equation}\label{eq:assocaction}
\xymatrix{
      x \ar[d] \ar[rrr]^{h \circ (g \circ f)}="0" \ar@/^3pc/[rrr]^{(h \circ g) \circ f}="1" & & & y \ar[d] \\
      x' \ar[rrr]_{h' \circ (g' \circ f')}  & \uriicell{(F \otimes_H G) \otimes_H H} & & y' \\
      \ar@{=>}"0" ;"1"^{a_{f,g,h}}
} \qquad = \qquad
\xymatrix{
      x \ar[rrr]^{(h \circ g) \circ f} \ar[d] & & & y \ar[d] \\
      x' \ar[rrr]_{(h' \circ g') \circ f}="0" \ar@/_3pc/[rrr]_{h' \circ (g' \circ f)}="1"  & \uriicell{F \otimes_H (G \otimes_H H)} & & y'\\
      \ar@{=>}"0" ;"1"^{a_{f',g',h'}}
}
\end{equation}
and similarly for vertical composition.  Likewise, unitors in the
horizontal and vertical bicategories agree with the identity for
composition of squares:
\begin{equation}\label{eq:unitaction}
\xymatrix{
      x \ar[d]_{g} \ar[r]_{f} \ar@/^2pc/[rr]^{f}="1" \ar@{}[rr]^{}="0" & y \ar[d]^{g'} \ar[r]_{1_y} & y \ar[d] \ar[d]^{g'} \\
      x' \ar[r]_{f'} \uriicell{F}  & y' \ar[r]_{1_{y'}} \uriicell{1_{g'}} & y' \\
      \ar@{=>}"0"+<0ex,+1.5ex>;"1"+<0ex,-1.5ex>^{l_f}
} \qquad = \qquad
\xymatrix{
      x \ar[d]_{g} \ar[r]^{1_x} &  \ar[d]_{g} \ar[r]^{f} & y \ar[d]^{g'} \\
      x' \ar[r]^{1_{x'}} \ar@/_2pc/[rr]_{f'}="1" \ar@{}[rr]^{}="0" \uriicell{1_{g}}  & x' \uriicell{F} \ar[r]^{f'} & y' \\
      \ar@{=>}"0"+<0ex,-1.5ex> ;"1"+<0ex,+2.5ex>^{r_{f'}}
}
\end{equation}
and similarly for vertical unitors.
\end{definition}

\begin{remark}This is rather more unwieldy than the definition of either a
bicategory or a double category, but is simpler than a
similarly elementary description of a {\db} (in the sense of Section
\ref{sec:internal}) would be.  In particular, where there are compatibility
conditions involving equations in this definition, a {\db} would have
only higher isomorphisms, themselves satisfying additional coherence
laws.  In particular, in {\vdbs}, the action of 2-morphisms on squares
is described by strict equations, rather than by a specified
isomorphism satisfying coherence laws.
\end{remark}

To help make sense of this definition, we note that it is possible
(following \cite[sec. 1.4]{verity}) to define categories
$\catname{Cyl_H}$ (respectively, $\catname{Cyl_V}$) of
\textit{cylinders}.  The objects of these categories are squares, and
maps are pairs of vertical (respectively, horizontal) 2-morphisms
joining the vertical (respectively, horizontal) source and targets of pairs of
squares which share the other two sides (this is shown in Table
\ref{table:doublebicatdata}, in Section \ref{sec:decatfy}: the
cylinders are ``thin'' versions of higher morphisms appearing there).
These are categories in the usual sense, with strict associativity and
unit laws.  These conditions would be weakened in a {\db} (in which
maps would include not just pairs of 2-morphisms, but also a
3-dimensional interior of the cylinder, which is a morphism in
$2\Mor$, or 2-morphism in $\Mor$, satisfying properties only up to a
4-dimensional 2-morphism in $2\Mor$).

\subsection{Decategorification}\label{sec:decatfy}

The main idea we are pursuing in this section is that of
``decategorification''.  This rather vague term refers to a process in
which category-theoretic information is discarded from a structure.
Typically, it refers to replacing isomorphisms with equations---for
example, a decategorification of the category of finite sets is the
set of their cardinalities, $\mathbbm{N}$.  Similarly, one can turn a
bicategory into a category by taking new morphisms to be 2-isomorphism
classes of old morphisms, and discarding all 2-morphisms.

To establish a relationship between the apparently conflicting terms
for the two types of ``{\db}'', we will now show how a {\vdb} can
arise as a decategorification of a {\db} satisfying some conditions.
We will show later that this can be done with the double cospan
examples of Section \ref{sec:dblspan}.  The conditions which are
needed allow us to speak of the ``action of 2-cells upon squares''.
To see what these are, we first consider a ``lower dimensional''
example of a similar process.  What we want to do to obtain {\vdbs}
has an analog in the case of double categories.

\begin{table}[h]
\begin{tabular}{|l|l|l|}
\hline
& $\Ob$ & $\M$ \\
\hline
Objects
&
\begin{minipage}{1in}
  \begin{equation*}
  \xymatrix{
  \bullet^{x} \\
  }
  \end{equation*}
\end{minipage}
 & 
\begin{minipage}{1in}
  \begin{equation*}
  \xymatrix{
  \bullet \ar[r]^{f} & \bullet
  }
  \end{equation*}
\end{minipage}
\\

\hline
Morphisms
&
\begin{minipage}{1in}

  \begin{equation*}
  \xymatrix{
  \bullet \ar[d]^{g} \\
  \bullet
  }
  \end{equation*}
\end{minipage}

 & 
\begin{minipage}{1in}
  \begin{equation*}
  \xymatrix{
  \bullet \ar[r] \ar[d] & \bullet \ar[d] \\
  \bullet \ar[r] \uriicell{F} & \bullet
  }
  \end{equation*}
\end{minipage}
\\

  \hline

\end{tabular}
\caption{Data of a Double Category\label{table:doublecatdata}}
\end{table}

In a double category, thought of as an internal category in $\Cat$, we
have data of four sorts, as shown in Table \ref{table:doublecatdata}.
That is, a double category $\catname{DC}$ has categories $\Ob$ of
objects and $\M$ of morphisms.  The first column of the table shows
the data of $\Ob$: its objects are the objects of $\catname{DC}$; its
morphisms are the \textit{vertical} morphisms.  The second column
shows the data of $\M$: its objects are the \textit{horizontal}
morphisms of $\catname{DC}$; its morphisms are the squares of
$\catname{DC}$.

There is a condition we can impose which
effectively turns the double category into a category, where the
horizontal and vertical morphisms are composable, and the squares can
be ignored.  The sort of condition involved is similar to the
\textit{horn-filling conditions} introduced by Ross Street
\cite{street} in his first introduction of the idea of weak
$\omega$-categories, or \textit{quasicategories}, in which all morphisms 
are $n$-simplexes for some $n$.  A horn filling condition says that,
given some hollow simplex with just one face (morphism) missing from the
boundary, there  will be a morphism to fill that face, and a compatible
``filler'' for the inside of the simplex.  In a double category, there is
an analogous ``niche-filler'' condition.

\begin{definition}A double category $\mathcal{DC}$ satisfies the
\textbf{composability condition} if the following holds.  For any pair
$(f,g)$ of a horizontal and vertical morphism where the target object of $f$
is the source object of $g$, there is a unique pair $(h,\star)$ consisting of
a unique vertical morphism $h$ and unique invertible square $\star$ making
the following diagram commute:
\begin{equation}\label{eq:hornfiller}
\xymatrix{
  x \ar@{-->}[d]_{h} \ar[r]^{f} & y \ar[d]^g \\
  z \ar[r]_{1_z} \uriicell{\star} & z
}
\end{equation}
and similarly when the source of $f$ is the target of $g$.
\end{definition}
Notice that taking $f$ to be the identity in this condition
implies $\star$ is the identity square.  This defines a composition:

\begin{theorem}\label{prop:fillercomp}If $\mathcal{DC}$ satisfies the
composability condition and there are no other squares in
$\mathcal{DC}$, there is a category $\mathcal{DC}_0$ with the same
objects as $\mathcal{DC}$, and all horizontal and vertical morphisms
as its morphisms.
\end{theorem}
\begin{proof}
Begin by defining composition from $\star$, so that if $f$ is
horizontal and $g$ is vertical, then $g\circ f = h$.  If $f$ and $g$
are both horizontal (or both vertical), then define $g \circ f$ to be
the usual composite.  Then this composition is associative and has
identities.  We only need to check this for the composition using
$\star$.  For example, given morphisms as in the diagram:
\begin{equation}
\xymatrix{
  w \ar[r]^{f} & x \ar[r]^{f'} & y \ar[d]^{g}\\
  z \ar[r]_{1_z} & z \ar[r]_{1_z} & z \\
}
\end{equation}
there are two ways to use the unique-filler principle to fill this
rectangle.  One way is to first compose the pairs of horizontal
morphisms on the top and bottom, then fill the resulting square.  The
square we get is unique, and the morphism is denoted $g \circ (f'
\circ f)$.  The second way is to first fill the right-hand square, and
then using the unique morphism we call $g \circ f'$, we get another
square on the left hand side, which our principle allows us to fill as
well.  The square is unique, and the resulting morphism is called $(g
\circ f') \circ f$.  Composing the two squares obtained this way must
give the square obtained the other way, since both make the diagram
commute, and both are unique.  So we have:
\begin{equation}\label{eq:fillercomp}
\xymatrix{
  w \ar@{-->}[d]_{(g \circ f') \circ f} \ar[r]^{f} & x \ar@{-->}[d] \ar[r]^{f'} & y \ar[d]^{g} \\
  z \ar[r]_{1_z} \uriicell{\star} & z \ar[r]_{1_z} \uriicell{\star}  & z \\
} \qquad = \qquad
\xymatrix{
  w \ar@{-->}[d]_{g \circ (f'\circ f)} \ar[r]^{f' \circ f} & y \ar[d]^{g} \\
  z \ar[r]_{1_z} \uriicell{\star} & z \\
}
\end{equation}
\end{proof}

\begin{remark} Note that the composability condition does not require
a square for every possible combination of source and target
morphisms.  In particular, there must be an identity morphism on the
boundary of the square---on the bottom in (\ref{eq:hornfiller}).  If
instead of the identity $1_z$, one could have any morphism $h$, then
by choosing $f$ and $g$ to be identities, this would imply that every
morphism must be invertible (at least weakly), since there must then
be an $h^{-1}$ with $h^{-1} \circ h$ isomorphic to the identity, but
of course we do not insist that all morphisms should have inverses.
When a filler square does exist, it indicates there is a commuting
square in $\mathcal{DC}_0$: the square $\star$ becomes an equation
between the composites along the upper right and lower left.
\end{remark}

The decategorification of a {\db} to give a {\vdb} is similar, except
that with a double category we were removing only the squares (the
lower-right quadrant of Table \ref{table:doublecatdata}).  There will
be a similar condition to satisfy, but we need to do more with a
{\db}, since there are more sorts of data (and, therefore, a more
complex condition).  These fall into a similar arrangement, as shown
in Table \ref{table:doublebicatdata}.

\begin{table}[h]
\begin{tabular}{|l|l|l|l|}
\hline
& $\Ob$ & $\M$ & $\B$ \\
\hline
Objects
&
\begin{minipage}{1in}
  \begin{equation*}
  \xymatrix{
  \bullet^{x} \\
  }
  \end{equation*}
\end{minipage}
 & 
\begin{minipage}{1in}
  \begin{equation*}
  \xymatrix{
  \bullet \ar[r]^{f} & \bullet
  }
  \end{equation*}
\end{minipage}
 & 
\begin{minipage}{1.5in}
  \begin{equation*}
  \xymatrix{
  \bullet \ar@/^1pc/[r]^{ }="0" \ar@/_1pc/[r]^{ }="1" & \bullet \\
  \ar@{=>}"0" ;"1"^{\alpha}
  }
  \end{equation*}
\end{minipage}
\\

\hline
Morphisms
&
\begin{minipage}{1in}

  \begin{equation*}
  \xymatrix{
  \bullet \ar[d]^{g} \\
  \bullet
  }
  \end{equation*}
\end{minipage}

& 
\begin{minipage}{1in}
  \begin{equation*}
  \xymatrix{
  \bullet \ar[r] \ar[d] & \bullet \ar[d] \\
  \bullet \ar[r] \uriicell{F} & \bullet
  }
  \end{equation*}
\end{minipage}
&
\begin{minipage}{1.5in}
  \begin{equation*}
  \xymatrix{
  \bullet \ar[d] \ar@/^1pc/[r]^{ }="0" \ar@/_1pc/[r]^{ }="1" & \bullet \ar[d] \\
  \ar@{=>}"0" ;"1"
  \bullet \uriicell{P_1} \ar@{-->}@/^1pc/[r]^{ }="2" \ar@/_1pc/[r]^{ }="3" & \bullet  \\
  \ar@{==>}"2" ;"3"
  }
  \end{equation*}
\end{minipage}

\\

  \hline
2-Cells
&
\begin{minipage}{1in}
  \begin{equation*}
  \xymatrix{
  \bullet \ar@/^1pc/[d]^{ }="0" \ar@/_1pc/[d]^{ }="1" \\
  \bullet \\
  \ar@{=>}"1" ;"0"^{\alpha}
  }
  \end{equation*}
\end{minipage}

&
\begin{minipage}{1in}
  \begin{equation*}
  \xymatrix{
  \bullet \ar@/^1pc/[d]^{ }="0" \ar@/_1pc/[d]^{ }="1" \ar[r] & \bullet  \ar@/^1pc/[d]^{ }="2" \ar@{-->}@/_1pc/[d]^{ }="3" \\
  \bullet \uriicell{P_2} \ar[r] & \bullet \\
  \ar@{=>}"1" ;"0"
  \ar@{==>}"3" ;"2"
  }
  \end{equation*}
\end{minipage}
&
\begin{minipage}{1.5in}
  \begin{equation*}
  \xymatrix{
  \bullet \ar@/^1pc/[dd]^{ }="0" \ar@/_1pc/[dd]^{ }="1" \ar@/^1pc/[rr]^{ }="2" \ar@/_1pc/[rr]^{ }="3" & & \bullet \ar@/^1pc/[dd]^{ }="4" \ar@{-->}@/_1pc/[dd]^{ }="5" \\
   & \Rightarrow^{T} & \\
  \bullet \ar@{-->}@/^1pc/[rr]^{ }="6" \ar@/_1pc/[rr]^{ }="7" & &
  \bullet \\ \ar@{=>}"1" ;"0" \ar@{=>}"2" ;"3" \ar@{==>}"5" ;"4"
  \ar@{==>}"6" ;"7" } \end{equation*}
\end{minipage}
  \\
\hline
\end{tabular}
\caption{The data of a \db\label{table:doublebicatdata}}
\end{table}

This table shows the data of the bicategories $\Ob$, $\M$, and $\B$,
each of which has objects, morphisms, and 2-cells.  Note that the
morphisms in the three entries in the lower right hand
corner---2-cells in $\M$, and morphisms and 2-cells in $\B$---are not
2-dimensional.  The 2-cells in $\M$ and morphisms in $\B$ are the
three-dimensional ``filling'' inside the illustrated cylinders, which
each have two square faces and two bigonal faces.  The 2-cells in $\B$
should be drawn as 4-dimensional.  The picture illustrated can be
thought of as taking both square faces of one cylinder $P_1$ to those
of another, $P_2$, by means of two other cylinders ($S_1$ and $S_2$,
say), in such a way that $P_1$ and $P_2$ share their bigonal faces.
This description works whether we consider the $P_i$ to be horizontal
and the $S_j$ vertical, or vice versa.  These describe the ``frame''
of this sort of morphism: the filling is the 4-dimensional ``track''
taking $P_1$ to $P_2$, or equivalently, $S_1$ to $S_2$, just as a
square in a double category can be read horizontally or vertically.
(Not all relevant parts of the diagrams have been labeled here, for
clarity.)

Next we want to describe a condition similar to the composability
condition for a double category.  In that case, we got a condition
which effectively allowed us to treat any square as an identity, so
that we only had objects and morphisms. Here, we want a condition
which lets us throw away the three entries of dimension greater than
two in Table \ref{table:doublebicatdata} in the bottom right.  This
condition, when satisfied, should allow us to treat a {\db} as a
{\vdb}.  It comes in three parts, one for each type of data we want to
discard:

\begin{definition}\label{def:actionconds} We say that a {\db} satisfies
the \textbf{vertical action condition} if, for any morphism $F_1 \in
\M$ and 2-morphism $\alpha \in \Ob$ such that $s(F_1) = t(\alpha)$,
there is a unique morphism $F_2 \in \M$ and unique invertible
2-morphism $P \in \M$ such that $P$ fills the ``pillow diagram'':
\begin{equation}\label{eq:vertactcond}
  \xymatrix{
    x \ar[r]^{ }="1" \ar[d] \ar@/^2pc/[r]^{ }="0" & y \ar[d] \\
    x' \ar[r] \uriicell{F_1} & y' \\
    \ar@{=>}"0" ;"1"^{\alpha}
  } \qquad \Rightarrow_P \qquad
  \xymatrix{
      x \ar [r] \ar[d] & y \ar[d] \\
      x' \ar[r]^{ }="3" \ar@/_2pc/[r]^{ }="2" \uriicell{F_2} & y' \\
      \ar@{=>}"3" ;"2"^{id}
  }
\end{equation} where $F_2$ is the back face of this diagram, and the
2-morphism in $\Ob$ at the bottom is the identity.

A {\db} satisfies the \textbf{horizontal action condition} if for any
morphism $F_1 \in \M$ and object $\alpha$ in $\B$ with $s(F_1) =
t(\alpha)$ there is a unique morphism $F_2 \in \M$ and unique
invertible morphism $P \in \B$ such that $P$ fill the pillow diagram:
\begin{equation}\label{eq:horizactcond}
  \xymatrix{
    x \ar[r] \ar[d]^{ }="1" \ar@/_2pc/[d]^{ }="0" & y \ar[d] \\
    x' \ar[r] \uriicell{F_1} & y' \\
    \ar@{=>}"0" ;"1"^{\alpha}
  } \qquad \Rightarrow_P \qquad
  \xymatrix{
      x \ar [r] \ar[d] & y \ar[d]^{ }="3" \ar@/^2pc/[d]^{ }="2"  \\
      x' \ar[r] \uriicell{F_2} & y' \\
      \ar@{=>}"3" ;"2"^{id}
  }
\end{equation}
\end{definition}

In (\ref{eq:vertactcond}), $F_2$ is the square which will eventually
be named $F_1 \star_H \alpha$ when we define an action of 2-cells on
squares, and in (\ref{eq:horizactcond}), $F_2$ is the square will
eventually be named $F_1 \star_V \alpha$.

\begin{remark} One can see that this condition is analogous to the
filler condition (\ref{eq:hornfiller}) in a double category by
imagining the diagram (\ref{eq:vertactcond}) viewed obliquely.  The
diagram says that given a square with two bigons---the top one
arbitrary and the bottom one the identity---there is another square
$F_2$ (the back face of a pillow diagram) and a filler 2-morphism $P
\in \B$ which fills the diagram.  If one imagines turning this diagram
on its side and viewing it obliquely, one sees precisely
(\ref{eq:hornfiller}), as a dimension has been suppressed.  The role
played by cylinders (2-morphism in $\B$) in (\ref{eq:vertactcond}) and
(\ref{eq:horizactcond} is played by a square in (\ref{eq:hornfiller});
the roles of both squares and bigons in (\ref{eq:vertactcond}) and
(\ref{eq:horizactcond}) are played by arrows in (\ref{eq:hornfiller});
the role of arrows in (\ref{eq:vertactcond}) and
(\ref{eq:horizactcond}) is filled by point-like objects in
(\ref{eq:hornfiller}).
\end{remark}

This gives horizontal and vertical actions, but to get the compatibility
between them, we need a further condition.  In particular, since these conditions 
involve both horizontal and vertical cylinders, the compatibility condition
must correspond to the 4-dimensional 2-cells in $\B$, shown in the lower
right corner of Table \ref{table:doublebicatdata}.

To draw the necessary condition is difficult, since the necessary
diagram is four-dimensional, but we can describe it as follows:

\begin{definition} We say a {\db} satisfies the
\textbf{action compatibility condition} if the following holds.
Suppose we are given \begin{itemize}
\item a morphism $F \in \M$
\item an object $\alpha \in \B$ whose target in $\M$ is a source of $F$
\item a 2-cell $\beta \in \Ob$ whose target morphism is a source of $F$
\item an invertible morphism $P_1 \in \B$ with $F$ as source, and the
      objects $\alpha$ and $\opname{id}$ in $\B$ as source and target
\item an invertible 2-cell $P_2 \in \M$ with $F$ as source, and the
      2-cells $\beta$ and $\opname{id}$ in $\M$ as source and target
\end{itemize} where $P_1$ and $P_2$ have, as targets, morphisms in
$\M$ we call $\alpha \star F$ and $\beta \star F$ respectively.  Then
there is a unique morphism $\hat{F}$ in $\M$ and unique invertible
2-cell $T$ in $\B$ having all of the above as sources and targets.
\end{definition}

Geometrically, the unique 2-cell in $\B$ looks like the structure in
the bottom right corner of Table \ref{table:doublebicatdata}.  This can be
seen as taking one horizontal cylinder to another in a way that fixes the (vertical)
bigons on its sides.  It does this by means of a translation which acts on the front
and back faces with a pair of vertical cylinders (have the same top
and bottom bigonal faces).  Alternatively, it can be seen as taking
one vertical cylinder to another, acting on the faces with a pair of
horizontal cylinders.  In either case, the cylinders involved in the
translation act on the faces, but the four-dimensional interior, $T$,
acts on the original cylinder to give another.  The simplest
interpretation of this condition is that it is precisely the condition
needed to give the compatibility condition (\ref{eq:actionindep}).

\begin{remark}Notice that the two conditions given imply the existence
of unique data of three different sorts in our \db.  If these are the
only data of these kinds, we can effectively omit them (since it
suffices to know information about their sources and targets).  This
omission is part of a decategorification of the same kind we saw for
a double category $\catname{DC}$.
\end{remark}

In particular, we show how a {\db} $\catname{D}$ satisfying the above
conditions gives a {\vdb}.  We know that $\catname{D}$ consists of
bicategories $(\Ob,\M,\B)$ together with all required maps (three
kinds of source and target maps, two kinds of identity, three
partially-defined compositions, left and right unitors, and the
associator), satisfying the usual properties.  To begin with, we
describe how the elements of a {\vdb} $\catname{V}$ (Definition
\ref{def:doublebicat}) arise from this:

\begin{definition}If $\catname{D}$ is a {\db}
satisfying the horizontal and vertical action conditions and the action
compatibility condition, then $\catname{V(D)}$ is the {\vdb} with:
\begin{itemize}
\item The objects $\Obj$ are the objects of $\Ob$.
\item The horizontal bicategory $\Hor$ of $\catname{V(D)}$ is $\Ob$
\item The vertical bicategory $\Ver$ of $\catname{V(D)}$ has:
  \begin{itemize}
  \item Objects: Objects of $\Ob$
  \item Morphisms: Objects of $\M$
  \item 2-morphisms: Objects of $\B$
  \end{itemize} The source, target and composition maps for $\Ver$ are
  the object maps from the source, target, and composition 2-functors
  for $\catname{D}$.
\item The squares $\Squ$ of $\catname{V(D)}$ are isomorphism classes
      of morphisms of $\M$.  These are equipped with:
  \begin{itemize}
  \item Vertical source and target maps: the morphism maps from the functors $s,t : \M \rightarrow \Ob$.
  \item Horizontal source and target maps: the internal ones in $\M$.
  \item Horizontal composition (\ref{eq:squarehorizcomp}): the composition of morphisms in $\M$.
  \item Vertical composition (\ref{eq:squarevertcomp}): the morphism maps for the partially defined functor $\circ$ for $\M$, 
  \item Horizontal Identity: The identity square for a morphism $g$ in
        $\Ver$ (i.e. $g$ an object in $\M$ is $1_f \in \M$.
  \item Vertical Identity: The identity square for a morphism $f$ in
        $\Hor$ (i.e. a morphism $f$ in $\Ob$) is given by
        $\opname{id}(f)$ for the unit functor $\opname{id} : \Ob
        \rightarrow \Mor$.
  \end{itemize}
\item The horizontal action defines $F \star_H \alpha$ to be (the
      isomorphism class of) the unique morphism in $\M$ whose
      existence is required by the horizontal action condition.
\item The vertical action defines $F \star_V \alpha$ to be (the
      isomorphism class of) the unique morphism in $\M$ whose
      existence is required by the vertical action condition.
\end{itemize}
\end{definition}

Of course, we must check this is really a {\vdb}:

\begin{theorem}\label{theorem:doublebicat}Suppose $\catname{D}$ is a {\db}
satisfying the horizontal and vertical action conditions and the action
compatibility condition.  Then $\catname{V(D)}$ is a {\vdb}.\end{theorem}
\begin{proof}
We check all the properties in the definition of a {\vdb}:
\begin{itemize}
\item By assumption, $\Hor$ is a bicategory.

\item $\Ver$ is a bicategory since the source and target functors in
      $\catname{D}$ for $\Ver$ satisfy all the usual axioms for a
      bicategory, hence their object maps do also.  Similarly, the
      composition maps have natural isomorphisms giving associators
      and unitors: they are just object maps of functors which satisfy
      the same conditions: in $\catname{D}$, the associator $a$
      satisfies the pentagon identity. The
      object maps for $a$ give the associator in $\Ver$.  Since the
      associator 2-natural transformation satisfies the pentagon
      identity, so do these object maps.  The other properties are
      shown similarly, so that $\Ver$ is a bicategory.

\item The source and target maps for $\Squ$ satisfies equations
      (\ref{eq:squarestmaps}) because the source and target maps of
      $\catname{D}$ are functors.

\item The composition lawsfor squares have the usual relation to
      source and target maps because, by assumption, $\M$ is a
      bicategory, but taking $\Squ$ to be 2-isomorphism classes of
      morphisms in $\M$, and disregarding all other 2-morphisms, we
      get that horizontal composition in $\Squ$ is exactly associative
      and has exact identities, so the squares are the morphisms of a
      category with respect to horizontal composition.

      Vertical composition for squares in $\catname{D}$ satisfies the
      axioms for a bicategory by the same argument as given above for
      $\Ver$, since it is the morphism map for the functor $\circ$. In
      particular, it has an associator and a unitor: but these must be
      morphisms in $\B$ since we take the morphism maps from the
      associator and unitor functors for $\circ$.  These are 2-isomorphisms,
      but since we defined squares to be 2-isomorphism classes (any
      isomorphism in $\B$ becomes an equation), this composition is exactly
      associative and has a unit.  Also, we disregard any morphisms in
      $\B$, so the squares are the morphisms of a category under
      vertical composition.

\item The interchange rule (\ref{eq:squareinterchangelaw}) follows
        from functoriality of the composition functors.

\item The actions $\star_H$ and $\star_V$ defined by the horizontal
      and vertical action conditions is well defined.  In particular,
      by composition of in $\M$ or $\B$, we guarantee the existence of
      the categories of horizontal and vertical cylinders
      $\catname{Cyl_H}$ and $\catname{Cyl_V}$, respectively.  These
      come from the 2-morphisms in $\M$ or morphisms in $\B$
      respectively which those conditions demand must exist.  Taking
      these to be identities, the cylinders consist of commuting
      cylindrical diagrams with two bigons and two squares.

      In the case where one bigon is the identity, and the other is
      any bigon $\alpha$, the conditions guarantee the existence of an
      invertible cylinder, which is now the identity because we have
      taken squares to be isomorphism classes.  This defines the
      effect of the action of $\alpha$ on the square whose source is
      the target of $\alpha$.  If this square is $F$, we denote the
      other square $\alpha \star_H F$ or $\alpha \star_V F$ as
      appropriate.

\item The horizontal action condition gives a well-defined action
      satisfies (\ref{eq:actioninterchange}) and
      (\ref{eq:actioncompat}) by an argument exactly analogous to that
      in the proof of Proposition \ref{prop:fillercomp}.  That is, the
      horizontal action condition means that certain fillers are
      unique.  When they can be obtained in two ways, these are equal.

\item The vertical action satisfies the vertical equivalent of
      (\ref{eq:actioninterchange}) and (\ref{eq:actioncompat}) for the
      same reason.

\item The condition (\ref{eq:actionindep}) guaranteeing independence
      of the horizontal and vertical actions follows from the action
      compatibility condition.  For suppose we have a square $F$ whose
      horizontal and vertical source arrows are the targets of 2-cells
      $\alpha$ and $\beta$, and attach to its opposite faces two
      identity 2-cells.  Then the horizontal and vertical action
      conditions mean that there will be a square $\alpha\star_H F$
      and a square $\beta\star_V F$).  Then the action compatibility
      condition applies (the $P_i$ are the identities we get from the
      action condition), and there is a morphism in $\M$, namely a
      square in $\catname{V}$ and a 2-cell $T \in \B$.  Consider the
      remaining face, which the action condition suggests we call
      $\alpha \star_H (\beta \star_V F)$ or $\beta \star_V (\alpha
      \star_H F)$, depending on the order in which we apply them.  The
      compatibility condition says that there is a unique square which
      fills this spot so the two must be equal.

\item We next check that composition for squares agrees with composition as in
      (\ref{eq:assocaction}).  Suppose we have three composable
      squares---that is, morphisms $F$, $G$, and $H$ in $\M$, which
      are composable along shared source and target objects in $\M$.
      The associator functor has an object map, giving objects in $\B$
      at the ``top'' and ``bottom'' of the squares.  It also has a
      morphism map, giving morphisms in $\B$.  But by assumption there
      is only a unique such map between , these associators must be
      the unique morphism in $\B$ with source $(H \circ G) \circ F$
      and target $H \circ (G \circ F)$.  Then by the vertical action
      condition, we have a filler 2-morphism in $\M$ for the action on
      the composite square by the top associator, and then, taking the
      result and composing with the bottom associator, we get another
      filler.  This must be the unique map between the two composites,
      which is the identity since they have the same sources and
      targets.  So we get a commuting cylinder.  Composing squares
      along source and target morphisms in $\Ob$ works the same way by
      a symmetric argument.

\item The condition (\ref{eq:unitaction}) is similar.  The unitor
      functor will give the unique morphism in $\B$, and the action
      compatibility condition gives the commuting cylinder for unitors
      on the composite of squares.

\end{itemize}

So indeed the construction of $\catname{V(D)}$ defines a {\vdb}.
\end{proof}

Next, in Section \ref{sec:equiv}, we continue the process of reducing
the complexity of these structures.  In particular, we see how {\vdbs}
can give rise to ordinary bicategories, which are frequently easier to
use.

\subsection{Bicategories from Double Bicategories}\label{sec:equiv}

It is well known that double categories can yield 2-categories in
three different ways.  Two obvious cases are when there are only
identity horizontal morphisms, or only identity vertical morphisms, so
that squares simply collapse into bigons with the two nontrivial
sides.  Notice that it is also true that a {\vdb} in which $\Hor$ is
trivial (equivalently, if $\Ver$ is trivial) is again a bicategory.
The squares become 2-morphisms in the obvious way, the action of
2-morphisms on squares is then just composition, and the composition
rules for squares in the double category become the rules for
composing 2-morphisms, and the result is clearly a bicategory.

The other, less obvious, case, is when the horizontal and vertical
categories on the objects are the same: this is the case of
\textit{path-symmetric} double categories, and the recovery of a
bicategory was shown by Brown and Spencer \cite{brownspencer}.  Fiore
\cite{fiore} shows how their demonstration of this fact is equivalent
to one involving \textit{folding structures}.

In this case we can interpret squares as bigons by composing the top
and right edges, and the left and bottom edges.  Introducing identity
bigons completes the structure.  These new bigons have a natural
composition inherited from that for squares.  It turns out that this
yields a bicategory.  Here, our goal will be to show half of an
analogous result, that a {\vdb} similarly gives rise to a bicategory
when the horizontal and vertical bicategories are equal.  We will also
show that a double bicategory for which the horizontal (or vertical)
bicategory is trivial can be seen as a bicategory.  The condition that
$\Hor = \Ver$ will hold in our general example of double cospans.

\begin{theorem}\label{thm:equiv} Any {\vdb} $\catname{V}=(\Obj, \Hor,
\Ver, \Squ, \otimes_H, \otimes_V, \star_H, \star_V)$ for which $\Hor =
\Ver$ produces a bicategory $\catname{B}$ by taking the 2-morphisms to be
2-morphisms in $\Hor$ and squares in $\Squ$.
\begin{proof}
We begin by defining the data of $\catname{B}$.  Its objects and
morphisms are the same as those of $\Hor$ (equivalently, $\Ver$).  We
describe the 2-morphisms by observing that $\catname{B}$ must contain
all those in $\Hor$ (equivalently, $\Ver$), but also some others,
which correspond to the squares in $\Squ$.

In particular, given a square
\begin{equation}\label{xy:sinsqu}
 \xymatrix{
  a \ar[r]^{f} \ar[d]_{g} & b \ar[d]^{g'} \\
  c \ar[r]_{f'} \uriicell{S} & d
 }
\end{equation}
there should be a 2-morphism
\begin{equation}
 \xymatrix{
   a \ar@/^1pc/[rr]^{g'\circ f}="0" \ar@/_1pc/[rr]_{f'\circ g}="1" && d \\
   \ar@{=>}"0"+<0ex,-2.5ex> ;"1"+<0ex,+2.5ex>^{S}
 }
\end{equation}

The composition of squares corresponds to either horizontal or
vertical composition of 2-morphisms in $\catname{B}$, and the relation
between these two is given in terms of the interchange law in a
bicategory:

Given a composite of squares,
\begin{equation}
  \xymatrix{
    x \ar[r]^{f} \ar[d]_{\phi_x} & y \ar[d]^{\phi_y} \ar[r]^{g} & z \ar[d]^{\phi_z} \\
    x' \ar[r]_{f'} \uriicell{F} & y' \ar[r]_{g'} \uriicell{G} & z'
 }
\end{equation}
there will be a corresponding diagram in $\catname{B}$:
\begin{equation}
 \xymatrix{
    x \ar@{}[rr]^{}="1"\ar[r]^{f} \ar@/_2pc/[rr]_{\phi_x \circ f'}="0" & y \ar@{}[rr]^{}="3" \ar[r]^{\phi_y} \ar@/^2pc/[rr]^{\phi_z \circ g}="2" & y' \ar[r]^{g'} & z' \\
   \ar@{=>}"1"+<0ex,-1.5ex> ;"0"+<0ex,+2.5ex>^{F}
   \ar@{=>}"2"+<0ex,-2.5ex> ;"3"+<0ex,+1.5ex>^{G}
 }
\end{equation}

Using horizontal composition with identity 2-morphisms
(``whiskering''), we can write this as a vertical composition:
\begin{equation}
 \xymatrix{
   x \ar@/^2pc/[rrr]^{\phi_z \circ g \circ f}="0" \ar[rrr]_{g' \circ \phi_y \circ f}="1" \ar@/_2pc/[rrr]_{g' \circ f' \circ \phi_x}="2" & & & z' \\
   \ar@{=>}"0"+<0ex,-2.5ex> ;"1"+<0ex,+2.5ex>^{G \circ \opname{1}_{f}}
   \ar@{=>}"1" ;"2"^{\opname{1}_{g'} \circ F}
 }
\end{equation}

So the square $F \otimes_H G$ corresponds to $(\opname{1} \circ G)
\cdot (F \circ \opname{1})$ for appropriate identities
$\opname{1}$.  Similarly, the vertical composite of $F' \otimes_V G'$
must be the same as $(\opname{1} \circ F) \cdot (G \circ
\opname{1})$.  Thus, every composite of squares which can be
built from horizontal and vertical composition, gives a corresponding
composite of 2-morphisms in $\catname{B}$, which are generated by
those corresponding to squares in $\Squ$, subject to the relations
imposed by the composition rules in a bicategory.

Now we want to show that {\vdb} $\catname{V}$ gives the entire
bicategory $\catname{B}$.  That is, that $\catname{B}$ has no other
2-morphisms than those which arise by the above process.  It suffices
to show that all such 2-morphisms not already in $\Hor$ arise as
squares (that is, the structure is closed under composition).  So
suppose we have any composable pair of 2-morphisms which arise from
squares $F$ and $G$.  If $F$ and $G$ have an edge in common, then we
have the situation depicted above (or possibly the corresponding form
in the vertical direction).  In this case, the composite 2-morphism
corresponds exactly to the composite of squares, and the axioms for
composition of squares ensure that all 2-morphisms generated this way
are already in our bicategory.  In particular, the unit squares become
unit 2-morphisms when composed with left and right unitors.

Now, if there is no edge in common to two squares, the 2-morphisms in
$\catname{B}$ must be made composable by whiskering with identities.
In this case, all the identities can be derived from 2-morphisms in
$\Hor$, or from identity squares in $\Squ$ (inside commuting
diagrams).  Clearly, any identity 2-morphism can be factored this way.
Then, again, the composite 2-morphisms in $\catname{B}$ will
correspond exactly to the composite of all such squares in $\Squ$ and
2-morphisms $\Hor$.

Finally, the associativity condition (\ref{eq:assocaction}) for the
action of 2-morphisms on squares ensures that composition of squares
agrees with that for 2-morphisms, so there are no extra squares from
composites of more than two squares.
\end{proof}
\end{theorem}
 
\begin{remark} When producing the bicategory $\catname{B}$ from
$\catname{V}$, we made a particular choice of orientation for the
2-morphisms obtained from squares.  The square $S \in \Squ$ shown in
(\ref{xy:sinsqu}) has vertical source $f$ and target $f'$, and
horizontal source $g$ and target $g'$.  However, the corresponding
2-morphism $S \in \catname{B}$ has source $g' \circ f$, which combines
vertical source and horizontal target; on the other hand, the target
of $S \in \catname{B}$ is $f' \circ g$, combining vertical target and
horizontal source.  We could equally well have chosen the opposite
convention.  This would give $\catname{B}^{\opname{co}}$, which is
$\catname{B}$ with the orientation of its 2-morphisms reversed.  (See,
e.g. \cite{leinster}).
\end{remark}

It is also worth considering here the situation of a double bicategory
in which all horizontal morphisms and 2-morphisms are identities.  In
this case, one can define a 2-morphism from a square with and bottom
edges being identities, whose source is the object whose identity is
the corresponding edge, and similarly for the target.  The composition
rules for squares in the vertical direction, then, are just the same
as those for a bicategory.  Likewise, the axioms for action of a
2-morphism on a square reduce to the composition laws for a bicategory
if one replaces the square by a 2-cell.

Next we describe a broad class of examples of {\dbs}, in the spirit of
the use of spans to give examples of bicategories.

\section{Double Cospans}\label{sec:dblspan}

In Remark \ref{thm:spanbicat} we described B\'enabou's demonstration
that $\Cspan$ is a bicategory for any category $\C$ with
pullbacks. Similarly, there is a bicategory of cospans in a category
$\C$ with pushouts.  There will be an analogous fact giving a {\db}
of double spans.  In fact, we describe this in terms of double
\textit{cospans}, since our aim in a subsequent paper will be to use
these to describe cobordisms, which have a natural description as
cospans.  Since cospans in $\C$ are the same as spans in the opposite
category, $\Cop$, this distinction is a matter of taste.

We remark here that similar constructions are described by Grandis
\cite{grandis3}, and related ``profunctor-based examples'' of
pseudo-double categories are described by Grandis and Par\'e
\cite{GP2}.  

\subsection{The Double Cospan Example}\label{sec:dcexample}

We begin by defining a {\db} of double cospans:

\begin{definition} $\iiCCosp$ is a {\db} of \textbf{double cospans} in $\C$,
consisting of the following: \begin{itemize}
\item the bicategory of objects is $\Ob = \CCosp$
\item the bicategory of morphisms $\M$ has: 
  \begin{itemize}
    \item as objects, cospans in $\C$;
    \item as morphisms, commuting diagrams of the form
          \begin{equation}\label{xy:cspan2}
            \xymatrix{
              X_1  \ar[r] \ar[d] & S \ar[d] & X_2 \ar[l] \ar[d] \\
              T_1 \ar[r]  & M  & T_2   \ar[l] \\
              X'_1 \ar[u] \ar[r] & S' \ar[u] & X'_2 \ar[u] \ar[l]
            }
          \end{equation}
          (in subsequent diagrams we suppress the labels for clarity);
    \item as 2-morphisms, cospans of cospan maps, namely commuting
      diagrams of the following shape:
        \begin{equation}\label{eq:cospanmap1}
        \xy
         (0,0)*{\bullet}="11";
         (0,15)*{\bullet}="12";
         (0,30)*{\bullet}="13";
         (15,0)*{\bullet}="21";
         (15,15)*{\bullet}="22";
         (15,30)*{\bullet}="23";
         (30,0)*{\bullet}="31";
         (30,15)*{\bullet}="32";
         (30,30)*{\bullet}="33";
         (-5,10)*{\bullet}="12a";
         (10,10)*{\bullet}="22a";
         (25,10)*{\bullet}="32a";
         {\ar "11";"12"};
         {\ar "13";"12"};
         {\ar "21";"22"};
         {\ar "23";"22"};
         {\ar "31";"32"};
         {\ar "33";"32"};
         {\ar "11";"21"};
         {\ar "31";"21"};
         {\ar "12";"22"};
         {\ar "32";"22"};
         {\ar "13";"23"};
         {\ar "33";"23"};
         {\ar "12a";"22a"};
         {\ar "32a";"22a"};
         {\ar "11";"12a"};
         {\ar "13";"12a"};
         {\ar "21";"22a"};
         {\ar "23";"22a"};
         {\ar "31";"32a"};
         {\ar "33";"32a"};
         {\ar "12";"12a"};
         {\ar "22";"22a"};
         {\ar "32";"32a"};
        \endxy
        \\
    \end{equation}
  \end{itemize}

\item the bicategory of 2-morphisms has:
  \begin{itemize}
  \item as objects, cospan maps in $\C$ as in (\ref{eq:spanmorph})
  \item as morphisms, cospan maps of cospans:
\begin{equation}\label{eq:cospanmap2}
\xy
 (0,0)*{\bullet}="11";
 (15,0)*{\bullet}="12";
 (30,0)*{\bullet}="13";
 (0,15)*{\bullet}="21";
 (15,15)*{\bullet}="22";
 (30,15)*{\bullet}="23";
 (0,30)*{\bullet}="31";
 (15,30)*{\bullet}="32";
 (30,30)*{\bullet}="33";
 (20,-5)*{\bullet}="12a";
 (20,10)*{\bullet}="22a";
 (20,25)*{\bullet}="32a";
     {\ar "11";"12"};
     {\ar "13";"12"};
     {\ar "21";"22"};
     {\ar "23";"22"};
     {\ar "31";"32"};
     {\ar "33";"32"};
     {\ar "11";"21"};
     {\ar "31";"21"};
     {\ar "12";"22"};
     {\ar "32";"22"};
     {\ar "13";"23"};
     {\ar "33";"23"};
     {\ar "12a";"22a"};
     {\ar "32a";"22a"};
     {\ar "11";"12a"};
     {\ar "13";"12a"};
     {\ar "21";"22a"};
     {\ar "23";"22a"};
     {\ar "31";"32a"};
     {\ar "33";"32a"};
     {\ar "12";"12a"};
     {\ar "22";"22a"};
     {\ar "32";"32a"};
\endxy
\\
\end{equation}

  \item as 2-morphisms, cospan maps of cospan maps:
\begin{equation}\label{eq:cospanmap3}
\xy
 (0,0)*{\bullet}="11";
 (30,0)*{\bullet}="12";
 (60,0)*{\bullet}="13";
 (0,30)*{\bullet}="21";
 (30,30)*{\bullet}="22";
 (60,30)*{\bullet}="23";
 (0,60)*{\bullet}="31";
 (30,60)*{\bullet}="32";
 (60,60)*{\bullet}="33";
 (40,-5)*{\bullet}="12a";
 (40,25)*{\bullet}="22a";
 (40,55)*{\bullet}="32a";
 (-5,20)*{\bullet}="21b";
 (25,20)*{\bullet}="22b";
 (55,20)*{\bullet}="23b";
 (35,15)*{\bullet}="22ab";
     {\ar "11";"12"};
     {\ar "13";"12"};
     {\ar "21";"22"};
     {\ar "23";"22"};
     {\ar "31";"32"};
     {\ar "33";"32"};
     {\ar "11";"21"};
     {\ar "31";"21"};
     {\ar "12";"22"};
     {\ar "32";"22"};
     {\ar "13";"23"};
     {\ar "33";"23"};
     {\ar "12a";"22a"};
     {\ar "32a";"22a"};
     {\ar "11";"12a"};
     {\ar "13";"12a"};
     {\ar "21";"22a"};
     {\ar "23";"22a"};
     {\ar "31";"32a"};
     {\ar "33";"32a"};
     {\ar "12";"12a"};
     {\ar "22";"22a"};
     {\ar "32";"32a"};
     {\ar "11";"21b"};
     {\ar "31";"21b"};
     {\ar "12";"22b"};
     {\ar "32";"22b"};
     {\ar "13";"23b"};
     {\ar "33";"23b"};
     {\ar "21b";"22ab"};
     {\ar "23b";"22ab"};
     {\ar "22b";"22ab"};
     {\ar "22a";"22ab"};
     {\ar "21";"21b"};
     {\ar "22";"22b"};
     {\ar "23";"23b"};
     {\ar "21b";"22b"};
     {\ar "23b";"22b"};
     {\ar "12a";"22ab"};
     {\ar "32a";"22ab"};
\endxy
\\
\end{equation}
  \end{itemize}
\end{itemize} All composition operations are by pushout; source and
target operations are the same as those for cospans.  The associators
and unitors in the horizontal and vertical bicategories are the maps
which come from the universal property of pushouts.
\end{definition}

\begin{remark} Just as 2-morphisms in $\M$ and morphisms in $\B$ can
be seen as diagrams which are ``products'' of a cospan with a map of
cospans, 2-morphisms in $\B$ are given by diagrams which are products
(as diagrams) of horizontal and vertical cospan maps.  These have, in
either direction, four maps of cospans, with objects joined by maps of
cospans.  Composition again is by pushout in composable pairs of
diagrams.
\end{remark}

Note that all these diagrams are products of smaller diagrams, each of
which is either a cospan, or a cospan map.  This suggests that the
horizontal and vertical directions should in some way behave like a
bicategory of cospans.  The next theorem shows this is indeed the case:

\begin{theorem}\label{theorem:span2bicat} For any category $\C$ with
pushouts, $\iiCCosp$ forms a {\db}.
\begin{proof}
$\M$ and $\B$ are bicategories since the composition functors act just
like composition in $\CCosp$, the bicategory of cospans in $\C$, in
each column, and therefore satisfies the same axioms.

Now, the horizontal and vertical directions have composition
operations defined in the same way.  Thus we can construct functors
between $\Ob$, $\M$, and $\B$ with the properties of a bicategory
simply by using the same constructions that turn each into a
bicategory in its own right.  In particular, the source and target
maps $s,t: \M \rightarrow \Ob$ and $s,t: \B \rightarrow \M$ are the
obvious maps giving the domains of the maps in (\ref{xy:cspan2}).  The
partially defined (horizontal) composition maps $\circ : \M^2
\rightarrow \M$ and $\otimes_H : \B^2\rightarrow \B$ are defined by
taking pushouts of diagrams in $\C$, which exist for any composable
pairs of diagrams because $\C$ has pushouts.  They are functorial
since they are independent of composition in the horizontal direction.
The associator for composition of morphisms is given in the pushout
construction.

To see that this construction gives a \db, we note that $\Ob$, $\M$,
and $\B$ as defined above are indeed bicategories.  Certainly, $\Ob$
is a bicategory because $\CCosp$ is a bicategory.  $\M$ and $\B$ are
bicategories because the morphism and 2-morphism maps from the
composition, associator, and other functors required for a {\db} give
them the structure of bicategories as well.

Moreover, the composition functors satisfy the properties of a
bicategory for just the same reason that composition of cospans (and
spans) does, since each of the three maps involved are given by this
construction.  Thus, we have a {\db}.
\end{proof}
\end{theorem}

Our motivation for Theorem \ref{theorem:span2bicat} is to show
that double cospans in suitable categories $\catname{C}$ give examples
of {\vdbs}.  We have described how to get a {\db} of such structures,
and we saw in Section \ref{sec:decatfy} that given certain conditions,
this gives a {\vdb}.  In Section \ref{sec:vdbcosp} we describe
explicitly the modifications we must make to $\CCosp$ to get these
conditions.

\subsection{A Verity Double Bicategory of Double Cospans}\label{sec:vdbcosp}

As we saw in Section \ref{sec:decatfy}, double bicategories have
higher morphisms of dimension up to $4$, but given certain conditions,
these can be omitted to give a {\vdb}.

\begin{definition}\label{def:cspan20}For a category $\C$ with
pushouts, the {\vdb} $\iiCCosp_0$, has: \begin{itemize}
\item the objects are objects of $\C$
\item the horizontal and vertical bicategories $\Hor = \Ver$ are both
      equal to a sub-bicategory of $\CCosp$, which includes only
      invertible cospan maps
\item the squares are isomorphism classes of commuting
      diagrams of the form (\ref{xy:cspan2})
\end{itemize} where two diagrams of the form (\ref{xy:cspan2}) are
isomorphic if they differ only in the middle objects, say $M$ and
$M'$, and the maps into these objects, and if there is an isomorphism $f:M
\rightarrow M'$ making the combined diagram commute.

The action of 2-morphisms $\alpha$ in $\Hor$ and $\Ver$ on squares is by
composition in diagrams of the form:
\begin{equation}\label{xy:cspan2action}
  \xymatrix{
      & \hat{S} \ar[d]^{\alpha} & \\
    X_1 \ar[ur]  \ar[r] \ar[d] & S \ar[d]^{s} & X_2 \ar[l] \ar[d] \ar[ul] \\
    T_1 \ar[r]^{t_1}  & M  & T_2 \ar[l]_{t_2} \\
    X'_1 \ar[u] \ar[r] & S' \ar[u]_{s'} & X'_2 \ar[u] \ar[l]
  }
\end{equation} (where the resulting square is as in \ref{xy:cspan2},
with $\hat{S}$ in place of $S$ and $s \circ \alpha$ in place of $s$).

Composition (horizontal or vertical) of squares of cospans is, as in
$\iiCCosp$, given by composition (by pushout) of the three cospans of
which the square is composed.  The composition for diagrams
of cospan maps are as usual in $\CCosp$.
\end{definition}

\begin{remark} Notice that $\Hor$ and $\Ver$ as defined are indeed
bicategories: eliminating all but the invertible 2-morphisms leaves a
collection which is closed under composition by pushouts.
\end{remark}

We will show more fully that this is a {\vdb} in Theorem
\ref{thm:cspanthm}.  First one must show that horizontal and vertical
composition of squares is well defined is defined on equivalence
classes.  We will get this result indirectly as a result of Theorems
\ref{theorem:span2bicat} and \ref{theorem:doublebicat}, but it is
instructive to see directly how this works in $\CCosp$.

\begin{theorem}In any category with pushouts, composition of squares
in Definition \ref{def:cspan20} is well-defined.
\begin{proof}
Suppose we have two representatives of a square, bounded by horizontal
cospans $X_1 \rightarrow S \leftarrow X_2$ and $X'_1 \rightarrow S'
\leftarrow X'_2$, and vertical cospans $X_1 \rightarrow T_1 \leftarrow X'_1$ 
and $X_2 \rightarrow T_2 \leftarrow X'_2$.  Suppose the middle objects
are $M$ and $\hat{M}$ as in the diagram (\ref{xy:cspan2}).  Given a
composable diagram which coincides along an edge (morphism in $\Hor$
or $\Ver$) with the first, we need to know that the two pushouts of
the different representatives are also isomorphic (that is, represent
the same composite square).

In the horizontal and vertical composition of these squares, the maps
to the middle object $M$ of the new square from the middle objects of
the new sides (given by composition of cospans) arise from the universal
property of the pushouts on the sides being composed (and the induced
maps from $M$ to the corners, via the maps in the cospans on the other
sides).  Since the middle objects are defined only up to isomorphism
class, so is the pushout: so the composition is well defined, since
the result is again a square of the form (\ref{xy:cspan2}).
\end{proof}
\end{theorem}

We use this, together with Theorems \ref{theorem:doublebicat} and
\ref{theorem:span2bicat}, (proved in Section \ref{sec:doublebicat}), to
show the following:

\begin{theorem}\label{thm:cspanthm} If $\C$ is a category with pushouts,
then $\iiCCosp_0$ is a {\vdb}.
\begin{proof}
In the construction of $\iiCCosp_0$, we take isomorphism classes of
double cospans as the squares.  We also restrict to invertible cospan maps
in the horizontal and vertical bicategories.

That is, take 2-isomorphism classes of morphisms in $\M$ in the \db,
where the 2-isomorphisms are invertible cospan maps, in both
horizontal and vertical directions.  We are then effectively
discarding all non-invertible morphisms and 2-morphisms in $\B$, and all
non-invertible 2-morphisms in $\M$.  In particular, there may be
``squares'' of the form (\ref{xy:cspan2}) in $\iiCCosp$ with
non-invertible maps joining their middle objects $M$, but we have
ignored these, and also ignore non-invertible cospan maps in the
horizontal and vertical bicategories.  Thus, we consider no diagrams of the form
(\ref{eq:cospanmap1}) except those in which the span maps are invertible,
in which case the middle objects are representatives of the same isomorphism
class. Similar reasoning applies to the 2-morphisms in $\B$.

The structure we get from discarding these will again be a
{\db}. In particular, the new $\M$ and $\B$ will be bicategories, since
they are, respectively, just a category and a set made into a discrete
bicategory by adding identity morphisms or 2-morphisms as needed.  On
the other hand, for the composition, source and target maps to be
bifunctors, the structures built from the objects, morphisms, and 2-cells
respectively must be bicategories.  This is since the composition, source,
and target maps are the object, morphism, and 2-morphism maps of these 
bifunctors, which satisfy the usual category axioms.  But the same argument
applies to those built from the morphisms and 2-cells as within $\M$ and
$\B$.  So we have a {\db}.

Next we show that the horizontal and vertical action conditions
(Definition \ref{def:actionconds} of Section \ref{sec:decatfy}) hold
in $\iiCCosp$.  A square in $\iiCCosp$ is a diagram of the form
(\ref{xy:cspan2}), and a 2-cell is a map of cospans.  Given a square
$M$ and 2-cell $\alpha$ with compatible source and targets as in the
action conditions, we have a diagram of the form shown in
(\ref{xy:cspan2action}).  Here, $M$ is the square diagram at the
bottom, whose top row is the cospan containing $S$.  The 2-cell
$\alpha$ is the cospan map including the arrow $\alpha: \hat{S}
\rightarrow S$.  There is a unique square built using the same objects
as $M$, but using the cospan containing $\hat{S}$ as the top row.  The
map from $\hat{S}$ to $M$ is then $s \circ \alpha$.

To satisfy the action condition, we want this square $\hat{M}$, which
is the candidate for $M_1 \star_V \alpha$, to be unique.  But suppose
there were another $\hat{M_2}$ with a map from $\hat{S}$.  Since we
are in $\iiCCosp_0$, $\alpha$ must be invertible, which would give a
map to $\hat{M_2}$ from $S$.  We then find that $\hat{M_2}$ and
$\hat{M}$ are representatives of the same isomorphism class, so in
fact this is the same square.  That is, there is a unique morphism in
$\B$ taking $M$ to $\hat{M}$ (a diagram of the form \ref{eq:cospanmap2})
with invertible cospan maps in the middle and bottom rows.  This is
the unique filler for the pillow diagram required by definition
\ref{def:actionconds}.

The argument that $\iiCCosp_0$ satisfies the action compatibility
condition is similar.

So $\iiCCosp_0$ is a {\db} in which, there there is at most one unique
morphism in $\M$, and at most unique morphisms and 2-morphisms in
$\B$, for any specified source and target, and the horizontal and
vertical action conditions hold.  So $\iiCCosp_0$ can be interpreted
as a {\vdb} (Theorem \ref{theorem:doublebicat}).
\end{proof}
\end{theorem}

\begin{remark}
We observe here that the compatibility condition (\ref{eq:assocaction})
relating the associator in the horizontal and vertical bicategories to
composition for squares is due to the fact that the associators
are maps which come from the universal property of pushouts.  This is
by the parallel argument to that we gave for spans in Remark
\ref{rem:spanbicat}.  The same argument applies to the middle objects
of the squares, and gives associator isomorphisms for that
composition.  When we reduce to isomorphism classes, these isomorphisms 
become identity maps, so we get a commuting pillow as in
(\ref{eq:assocaction}).  A similar argument shows the compatibility
condition for the unitor, (\ref{eq:unitaction}).
\end{remark}

It is interesting to note how the arguments in the proof of
Theorem \ref{theorem:doublebicat} apply to the case of $\iiCCosp$.

In particular, the interchange rules hold because the middle objects
in the four squares being composed form the vertices of a new square.
The pushouts in the vertical and horizontal direction form the
middle objects of vertical and horizontal cospans over these.  The
interchange law means that the pushout (in the horizontal direction)
of the objects from the vertical cospans is in the same isomorphism
class as the pushout (in the vertical direction) of the objects from
the horizontal cospans.  This follows from the universal property
of the pushout.

\subsection{Example: Cobordisms with Corners}

One important example of a category of cospans involves cobordism of
manifolds, although to realize this example requires some additional
structure.  In particular, the category $\nCob$ of cobordisms with
corners is not $\iiCCosp$ for a category $\C$ with pushouts, since the
objects of this category are manifolds, and $\catname{Man}$ does not
have pushouts.

Recall that two manifolds $S_1, S_2$ are \textit{cobordant} if there
is a compact manifold with boundary, $M$, such that $\partial M$ is
isomorphic to the disjoint union of $S_1$ and $S_2$.  This $M$ is
called a \textit{cobordism} between $S_1$ and $S_2$.  So in
particular, a cobordism is a cospan $S_1 \rightarrow M \leftarrow
S_2$, where both maps are inclusions of the boundary components.  A
cobordism with corners is then a manifold with corners, where the
boundary components and corners are included in a double cospan.

In particular, for topological cobordisms (i.e. cobordisms which are
topological manifolds with boundary), all the pushouts required to
compose such double cospans will still be topological manifolds.  For
smooth manifolds, to ensure that the result of gluing is smooth we
need to specify an additional condition, using "collars" on the
boundaries and corners.

\subsubsection{Cobordisms and Collars}\label{sec:collarman}

To begin with, recall that a smooth manifold with corners is a
topological manifold with boundary, together with a maximal compatible
set of coordinate charts $\phi : \Omega \rightarrow [0,\infty)^n$ -
into the positive sector of $latex \mathbb{R}^n$. (where $\phi_1$,
$\phi_2$ are compatible if $\phi_2 \circ \phi_1^{-1}$ is a
diffeomorphism).

J\"anich \cite{jan68} introduces the notion of $\br{n}$-manifold,
reviewed by Laures \cite{laures}.  This is build on a manifold with
corners, using the notion of a \textit{face}:

\begin{definition}(\textbf{J\"anich})A \textbf{face} of a manifold
with corners is the closure of some connected component of the set of
points with exactly one zero component in any coordinate chart.  An
$\br{n}$-\textit{manifold} is a manifold with faces together with an
$n$-tuple $(\partial_0 M, \dots, \partial_{n-1}M)$ of faces of $M$,
such that
\begin{itemize}
\item $\partial_0 M \cup \dots \partial_{n-1} M = \partial M$
\item $\partial_i M \cap \partial_j M$ is a face of $\partial_i M$ and
      $\partial_j M$
\end{itemize}
\end{definition}

The case we will be interested in here is the case of
$\br{2}$-manifolds.  In this notation, a $\br{0}$-manifold is just a
manifold without boundary, a $\br{1}$-manifold is a manifold with
boundary, and a $\br{2}$-manifold is a manifold with corners whose
boundary decomposes into two components (of codimension 1), whose
intersections form the corners (of codimension 2).  We can think of
$\partial_0 M$ and $\partial_1 M$ as the ``horizontal'' and
``vertical'' part of the boundary of $M$.

Now, for a point $x \in S$, there will be a neighborhood $U$
of $x$ which restricts to $U_1 \subset M_1$ and $U_2 \subset M_2$ with
smooth maps $\phi_i : U_i \rightarrow [0,\infty)^n$ with $\phi_i(x)$
on the boundary of $[0,\infty)^n$ with exactly one coordinate equal to
$0$.  One can easily combine these to give a homeomorphism $\phi: U
\rightarrow \mathbbm{R}^n$, but this will not necessarily be a
diffeomorphism along the boundary $S$.  While \textit{topological}
cobordisms can be composed along their boundaries, \textit{smooth}
cobordisms $M_1$ and $M_2$ should be composed differently, to ensure
that every point---including points on the boundary of $M_i$---will
have a neighborhood with a smooth coordinate chart.  To solve this
problem, we use \textit{collars}, which is also done in the category
$\nCobi$.

The \textit{collaring theorem} says that for any smooth manifold with
boundary $M$, $\partial M$ has a \textit{collar}: an embedding $f :
\partial M \times [0,\infty) \rightarrow M$, with $(x,0) \mapsto x$ for $x \in \partial M$.  
This is a well-known result (for a proof, see e.g. \cite{hirsch}, sec.
4.6).  It is an easy corollary that we can choose
to use the interval $[0,1]$ in place of $[0,\infty)$ here.

Laures (\cite{laures}, Lemma 2.1.6) describes a generalization of this
theorem to $\br{n}$-manifolds, so that for any $\br{n}$-manifold $M$,
there is an $n$-dimensional cubical diagram
($\br{n}$-\textit{diagram}) of embeddings of cornered neighborhoods of
the faces. 
It is then standard that one can compose two smooth cobordisms
with corners, equipped with such smooth collars, by gluing along $S$.
The composite is then the topological pushout of the two inclusions.
Along the collars of $S$ in $M_1$ and $M_2$, charts $\phi_i : U_i
\rightarrow [0,\infty)^n$ are equivalent to charts mapping into
$\mathbbm{R}^{n-1} \times [0,\infty)$, and since the the composite has
a smooth structure defined up to a diffeomorphism which is the
identity along $S$.  The precise smooth structure on this cobordism
depends on the collar chosen, but there is always such a choice, and
the resulting composites are all equivalent up to diffeomorphism.
 
Now, for each $n$, one can define:

\begin{definition}The bicategory $\nCob$ is given by the following data:
\begin{itemize}
\item The objects of $\nCob$ are of the form $P = \hat{P} \times I^2$
      where $\hat{P}$ may be any $(n-2)$ manifolds without boundary
      and $I=[0,1]$.
\item The morphisms of $\nCob$ are cobordisms $P_1
      \rightarrowlim^{i_1} S \leftarrowlim^{i_2} P_2$
        where $S = \hat{S} \times I$ and $\hat{S}$ is an
        $(n-1)$-dimensional collared cobordism with corners such that: the
        $\hat{P_i} \times I$ are objects, the maps are injections into
        $S$, a manifold with boundary, such that $i_1(P_1) \cup
        i_2(P_2) = \partial S \times I$, $i_1(P_1) \cap i_2(P_2) =
        \emptyset$,
\item The 2-morphisms of $\nCob$ are generated by:
  \begin{itemize}
    \item diffeomorphisms of the form $f \times \opname{id}: T\times[0,1] \rightarrow T' \times
        [0,1]$ where $T$ and $T'$ have a common boundary, and $f$ is a
        diffeomorphism$T \rightarrow T'$ compatible with the source
        and target maps, i.e. fixing the collar.
    \item 2-cells: diffeomorphism classes of $n$-dimensional manifolds $M$
      with corners satisfying the properties of $M$ in the diagram of
      equation (\ref{xy:cspan2}), where isomorphisms are
      diffeomorphisms preserving the boundary
  \end{itemize}
      where the composite of the diffeomorphisms with the 2-cells
      (classes of manifolds $M$) is given by composition of
      diffeomorphisms of the boundary cobordisms with the injection
      maps of the boundary $M$
\end{itemize} The source and target objects of any cobordism $S$ are
specified by saying that the source of $S$ is the collection of
components of $\partial S \times I$ for which the image of $(x,0)$
lies on the boundary for $x \in \partial S$, and the target has the
image of $(x,1)$.  The source and target objects are the collars,
embedded in the cobordism in such a way that the source object $P =
\hat{P} \times I^2$ is embedded in the cobordism $S =
\hat{S} \times I$ by a map which is the identity on $I$ taking the
first interval in the object to the interval for a horizontal
morphism, and the second to the interval for a vertical morphism.  The
same condition distinguishing source and target applies as above.

Composition of 2-cells works by gluing along common boundaries.
\end{definition}

\begin{lemma}\label{lemma:ncobclosedcomp} Composition of morphisms and
2-morphisms in $\nCob$ is well-defined and $\nCob$ is closed under
composition.
\begin{proof}
The horizontal and vertical morphisms are products of the interval $I$
with $\br{1}$-manifolds, whose boundary is $\partial_0 S$), equipped
with collars.  Suppose we are given two such cobordisms $S_1$ and
$S_2$, and an identification of the source of $S_2$ with the target of
$S_1$ (say this is $P = \hat{P} \times I$).  Then the composite $S_2
\circ S_1$ is topologically the pushout of $S_1$ and $S_2$ over $P$.
Now, $P$ is smoothly embedded in $S_1$ and $S_2$, and any point in the
pushout will be in the interior of either $S_1$ or $S_2$ since for any
point on $\hat{P}$ each end of the interval $I$ occurs as the boundary
of only one of the two cobordisms.  So the result is smooth.  Thus,
$\iiCob$ is closed under such composition of morphisms.

The same argument holds for 2-cells, since it holds for any
representative of the equivalence class of some manifold with corners,
$M$, and the differentiable structure will be the same, since we
consider equivalence up to diffeomorphisms which preserve the collar
exactly.
\end{proof}
\end{lemma}

\begin{figure}[h]
\begin{center}
\includegraphics{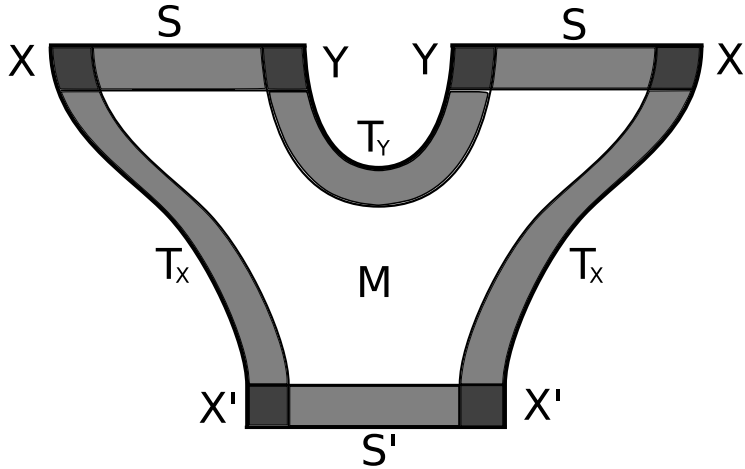}
\end{center}
\caption{\label{fig:cobcorners-labelled}A 2-Dimensional Cobordism with Corners}
\end{figure}

\begin{figure}[h]
\begin{center}
\includegraphics{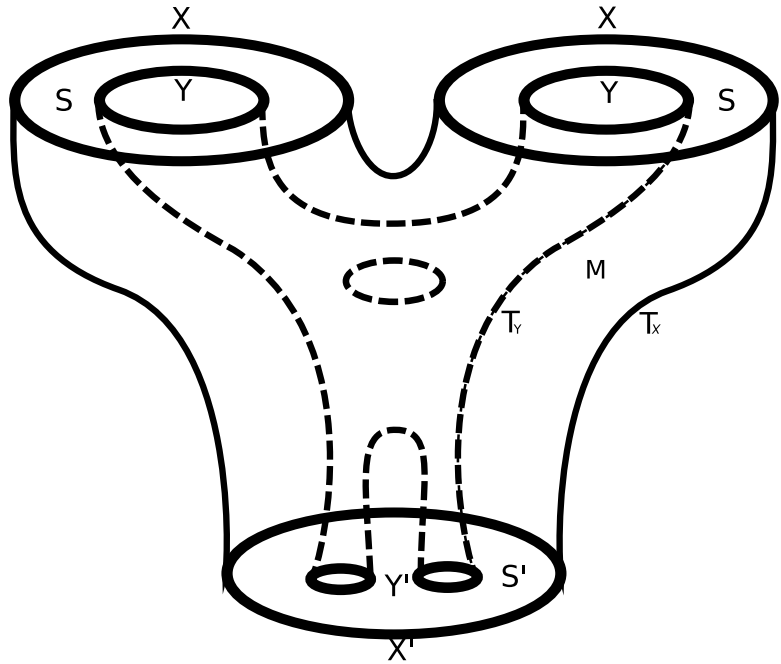}
\end{center}
\caption{\label{fig:cobcorn3d}A 3-Dimensional Cobordism with Corners}
\end{figure}

Examples of such cobordisms with corners in 2 and 3 dimensions, as
illustrations of (\ref{xy:cspan2}), are shown in in Figures
\ref{fig:cobcorners-labelled} and \ref{fig:cobcorn3d}, respectively.

\subsection{Prospects for $n$-tuple Bicategories}\label{sec:ntplbicat}

We have discussed both {\dbs} and {\vdbs}, both of which can be seen
as forms of \textit{weak cubical} $n$-category.  The broader question
of various definitions of weak $n$-category is discussed in more
detail by Tom Leinster \cite{leinster}, and by Eugenia Cheng and Aaron
Lauda \cite{chenglauda}.  In light of this context, Theorem
\ref{theorem:span2bicat} suggests one direction of generalization for
{\dbs}, to ``$n$-tuple bicategories'' for any $n$.  Moreover, our
example of double (co)spans can be generalized to arbitrarily high
dimension. 

We have seen how to construct $\iiCCosp$ for a category $\C$ with
pushouts, and how we take a restricted form of this construction to
yield a {\vdb}.  We have chosen to stop the process of
taking cospans in a category of cospans after two steps, but we could
continue this construction.  Taking cospans in this new category gives
cubes of objects with maps from corners to the middles of edges, from
middles of edges to middles of faces, and from middles of faces to the
middle of the cube.  Similarly, for any finite $n$, we can iterate the
process of taking cospans to yield an $n$-dimensional cube.

In particular, we note that ``{\vdbs}'' arise from special examples of
bicategories internal to $\Bicat$.  There is a category of all such
structures, namely the functor category of all maps $F: Th(\Bicat)
\rightarrow \Bicat$, denoted $[ Th(\Bicat), \Bicat ]$.  There will be
an analogous concept of ``triple bicategories'', namely bicategories
internal to $[ Th(\Bicat), \Bicat ]$.  In general, a ``$k$-tuple
bicategory'' will be a bicategory internal to the category of weak
$(k-1)$-tuple categories.

We expect that for all $k$, a $k$-tuply iterated process of taking
cospans of cospans (or similarly for spans) will yield examples of
these structures.  These $k$-dimensional (co)spans will naturally form
a weak $k$-tuple category.  Marco Grandis \cite{grandis3} describes
this in terms of a somewhat different description of weak $n$-cubical
categories.

A further direction of generalization would be to substitute
tricategories, tetracategories, and so forth in place of bicategories
in the preceding construction, perhaps making different choices each
stage.  The question then arises what sort of structures it would be
possible to define by selectively decategorifying, and what sorts of
``filler'' conditions this would need.  Another potentially
interesting question is whether the examples based on cospans also
generalize---perhaps by taking cospans, not in a category, but in an
$n$-category.

\section{Acknowledgments}

The author would like to acknowledge the help of John Baez in forming
this project.  Thanks and recognition are due to Dominic Verity for
providing, in his Ph.D. thesis, the crucial definition; to Marco Grandis
for significantly developing of the theory of cubical cospans since the
original version of this paper was written; to Aaron Lauda for references 
on manifolds with corners; and to Tom Fiore for helpful discussion on weak 
double categories.  Acknowledgment is also due to Dan Christensen, James Dolan, Derek Wise, Toby Bartels, Mike Stay, Alex Hoffnung, and John Huerta for discussions of the work in progress.

\end{document}